\documentclass{amsart}
\usepackage{graphicx}
\graphicspath{{Fig_file/}}
\usepackage{amsfonts}
\usepackage{amssymb}
\usepackage{mathrsfs}
\usepackage{hyperref}
\usepackage{xypic}
\usepackage{cite}
\usepackage{float}
\usepackage{subfigure}
\newtheorem{thm}{Theorem}[section]
\newtheorem{cor}[thm]{Corollary}
\newtheorem{lem}[thm]{Lemma}

\theoremstyle{definition}

\theoremstyle{remark}
\newtheorem{rem}[thm]{Remark}
\numberwithin{equation}{section}

\begin{document}
\title[Y.S. Jiang]{Modeling  Boyciana-fish-human Interaction with Partial Differential Algebraic Equations}
\author[Y.S. Jiang]{Yushan Jiang$^\ast$}
\author[Q.L. Zhang]{Qingling Zhang}
\author[H.Y. Wang]{Hanyan Wang}
\address[Y.S. Jiang, Q.L. Zhang]{Institute of System Science Northeastern University, Shenyang,China.}%
\address[H.Y. Wang]{School of Mathematical $\&$ Natural Sciences, Arizona State University, U.S.A.}%
\email{sobolev@126.com,jys@neuq.edu.cn}%
\begin{abstract}
  With human social behaviors influence, some boyciana-fish reaction-diffusion system coupled with elliptic human distribution equation is considered. Firstly, under homogeneous Neumann boundary conditions and ratio-dependent functional response the system can be described as a nonlinear partial differential algebraic equations (PDAEs) and the corresponding linearized system is discussed with singular system theorem. In what follows we discuss the elliptic subsystem and show that the three kinds of nonnegative are corresponded to three different human interference conditions: human free, overdevelopment and regular human activity. Next we examine the system persistence properties: absorbtion region and the stability of positive steady states of three systems. And the diffusion-driven unstable property is also discussed. Moreover, we propose some energy estimation discussion to reveal the dynamic property among the boyciana-fish-human interaction systems.Finally, using the realistic data collected in the past fourteen years, by PDAEs model parameter optimization, we carry out some predicted results about wetland boyciana population. The applicability of the proposed approaches are confirmed analytically and are evaluated in numerical simulations.
\end{abstract}
\maketitle
\section{Introduction}
In Northern China, Beidaihe wetland is located in the junction of three big ecosystems that are forest, ocean and wetland. Beidaihe wetland is one of the important channels for far east migratory birds. Every September, October, there are about 400 species of birds migrate to Beidaihe wetland park\cite{2015HongyingGao}. At the same time, as an international tourist city Beidaihe attracted a large number of tourists to travel from May to October. Since the increasing human activities, many living and breeding areas were damaged and lead to the decrease of the bird population.

 Boyciana is one the most sensitive species in Beidaihe wetland system. It is the first-grade state protection animal in China. Boyciana was widely distributed in Northeast Asia. In 1986, the experts found that 2729 boyciana moved through Beidaihe wetland region. However, in recent decades the human's activities made boyciana's predatory object quantity to reducing. The habitat environment of boyciana was also destroyed. Many environmental factors influence the spatiotemporal distribution of boyciana such as hidden factor, water factor, vegetation factor and food factor are directly or indirectly related to human activities\cite{2014XiaoxuWu}.

 Despite a rich literature on the spatiotemporal research of ecosystem\cite{2014Guichard}. The human's interference to the ecological spatiotemporal process is rarely studied. Thus, the goal of our theoretical ecology model is to study how the interactions between boyciana and the their food-wetland fish with human social behavior influence.

Mathematically, reaction-diffusion equation can be used to model the spatiotemporal distribution and abundance of organisms\cite{2015JiantaoZhao_WeiJJ,2015KengDeng,2015GuoweiDai_WangHY,2013ZhanpingMa}. In recent decades the role of the reaction-diffusion effect in maintaining bio-diversity has received a great deal of attention in the literature on ecology conservation\cite{2015Mimura}. Empirical evidence suggests that the spatial scale can influence population interactions.
The major class of spatial models are those that treat space as a continuum and describe the distribution of populations in terms of densities. A typical form of reaction-diffusion population model is
$$\frac{\partial x}{\partial t}=D\Delta x+xf(z, x)$$
where $x(t,z)$ is the population densities vector at time $t$ and space point $z$, $D$ is the diffusion constant matrix, $\Delta$ is the Laplace operator with respect to the spatial variable $z$, and $f(z, x)$ is the growth function vector.
Such an ecological model was first considered by Skellam \cite{1991Skellam}, and the reaction-diffusion biological models were also studied by Fisher \cite{1937FISHER} and Kolmogoroff\cite{1937Kolmogorov}.

For the predator-prey type reaction-diffusion biological models,  we are referred to the functional response models traditionally used as Lotka-Volterra \cite{2011Galiano}, Allee effect \cite{2015XuechenWang_WeiJJ}, Holloing type \cite{2011ZijianLiu}, Bedding-DeAngelis \cite{2014XiangpingYan}, ratio-dependent \cite{2013Ko} et al. All the predator-prey models cannot be directly applied in our human interference model (\ref{prey_predator systems}). Our developed system includes three interaction species: boyciana, fish and human.

To this end, we modify the ratio-dependent reaction-diffusion system model \cite{2013Ko} in which the incorporating one prey and two competing predator species was considered. In \cite{2013Ko} we replace one competing predator species by human. Specially, the human influence part is some degenerated Fisher population model of elliptic type. Although the global attractor and persistence of the system was discussed in \cite{2013Ko} by comparison principle theory for parabolic equations. While our model involves some elliptic equation, so we propose some new method on system persistence property. The diffusion-driven instability or turing instability which has attracted the attention of some investigators is also discussed in this study by using the qualitative theory.

\section{Boyciana-fish-human Model}
\subsection{Mathematical model}
  In this section, we will propose some singular PDEs system (or PDAEs system) on wetland ecological system as the following form:
$$E\frac{\partial x}{\partial t}=D\Delta x+xf(\textbf{z}, x)$$
where $E$ is a singular matrix, $D$ is the diffusion matrix and $f$ is the ratio dependent response functional function vector. We choose human, boyciana and fish as the research objects. The spatiotemporal dynamics between boyciana and fish with human activity affect in a protected environment can be described by the following PDAEs
\begin{equation}\label{prey_predator systems}
\left\{
\begin{aligned}
  &\frac{\partial x_1}{\partial t}
  =d_{1}\Delta x_1+x_1\left(1-x_1-\frac{cx_2}{x_1+\alpha x_2}- h_1x_3\right), \textbf{z}\in\Omega,t>0,\\
  &\frac{\partial x_2}{\partial t}
  =d_{2}\Delta x_2+x_2\left(-d+\frac{mx_1}{x_1+\alpha x_2}-h_2x_3\right), \textbf{z}\in\Omega,t>0,\\
  &0
  =\Delta x_3+rx_3(1-x_3), \textbf{z}\in\Omega,\\
  &\frac{\partial x_1}{\partial \textbf{n}}
  =\frac{\partial x_2}{\partial \textbf{n}}
  =\frac{\partial x_3}{\partial \textbf{n}}=0, \textbf{z}\in\partial\Omega,t>0,\\
     &x_1(0,\textbf{z})=x_1^0(\textbf{z})\geq 0,x_2(0,\textbf{z})=x_2^0(\textbf{z})\geq 0 , \textbf{z}\in\Omega.\\
\end{aligned}
\right.
\end{equation}
where $x_1(t,\textbf{z}), x_2(t,\textbf{z}), x_3(\textbf{z})$ are the state variables, $\Omega$ is a bounded spatial region, $\textbf{z}=[z_1, z_2]\in \Omega$ is the spatial coordinate, $\textbf{n} $ is the outward unit normal vector of the boundary $\partial\Omega$, the coefficients $c, \alpha, m, d, h_1, h_2$ and $r$  are positive constants. The initial value $x_1^0(\textbf{z}),x_2^0(\textbf{z})$
are non-negative smooth functions which are not identically zero.

By partial differential algebraic equations (PDAEs) theory \cite{Jiang2015}, we can rewrite (\ref{prey_predator systems}) as the following matrix form
\begin{equation}\label{PDAE01}
  E\frac{\partial {x}}{\partial t}=D\Delta{x}+{f}({x})
\end{equation}
subject to the boundary condition (BC)
    \begin{equation}\label{PDAE02}
    \frac{\partial x_\iota(t,\textbf{z})}{\partial \textbf{n}}=0, \textbf{z}\in\partial\Omega,\iota\in\vartheta_{BC}
    \end{equation}
and the initial condition (IC)
    \begin{equation}\label{PDAE03}
    x_\kappa(0,\textbf{z})=x^0_\kappa(\textbf{z}),\textbf{z}\in\Omega,\kappa\in\vartheta_{IC}
    \end{equation}
where $x=(x_1,x_2,x_3)^T$, $\vartheta_{BC}=\{1,2,3\}$, $\vartheta_{IC}=\{1,2\}$ and
$$E=\left(
           \begin{array}{ccc}
             1 & 0 & 0 \\
             0 & 1 & 0 \\
             0 & 0 & 0 \\
           \end{array}
         \right),
D=\textrm{diag}(D_1,1),D_1=\left(
        \begin{array}{cc}
          d_{1} & 0 \\
          0 & d_{2} \\
        \end{array}
      \right),
$$
$${f}({x})=\left(
                    \begin{array}{c}
                      x_1 \left( 1-x_1-\frac{cx_2}{x_1+\alpha x_2}- h_1x_3\right) \\
                      x_2\left(-b+\frac{mx_1}{x_1+\alpha x_2}-h_2x_3\right) \\
                      rx_3(1-x_3)
                    \end{array}
                  \right).
$$
 The applications and mathematically research on PDAEs have attracted increasing attention to academics\cite{2014HuainingWu,2014Daafouz,2011ShuxiaTang}. In view of the latest literature in reaction-diffusion system research, most of the derivative coefficient matrix $D$ of these system is invertible. In other words, they are the parabolic type nonlinear partial differential equations system (PDEs). Most of them can be analyzed with the existing qualitative theory \cite{2012Pao} directly. However, for some parabolic-elliptic system \cite{2014LiangchenWang} or singular system \cite{Zhang2012} with reaction-diffusion term, the existing research results are relative few. If without considering the effects of spatial variables, the above system becomes a singular system or generalized state-space system\cite{Yang2013}.
 From the above PDAEs system (\ref{PDAE01})-(\ref{PDAE03}) we know the time derivative coefficient matrix $E$ is singular. Mathematically, it is a generalization of the classical parabolic PDEs system. Because some theoretical results \cite{Smoller1994,1998Christofides} can not be direct application of the singular situation, research in the singular case is relatively scarce. Some identifiability and stability properties has been studied in \cite{Jiang2015} with singular system theory. In this study, we proposed some theoretical results on this PDAEs system.

\subsection{Ecological description}
\begin{figure}
  \includegraphics[width=0.8\textwidth]{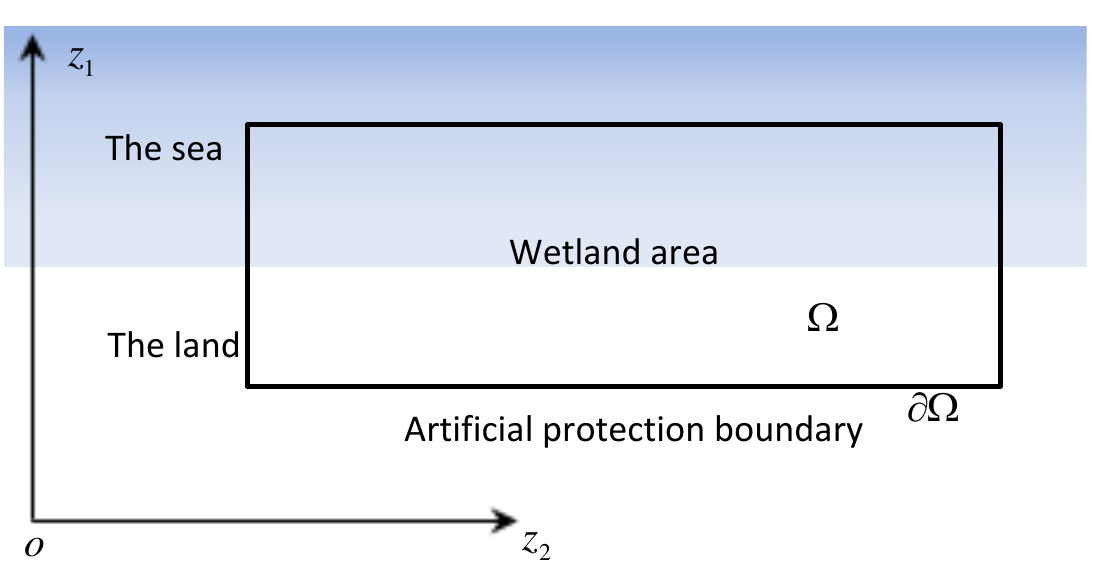}\\
  \caption{\label{ecologysys}The idealized spatial domain is a rectangular domain with the sea oriented direction $z_1$ and coast line direction $z_2$. The wetland conservation is closed with no flux boundary conditions imposed.}
\end{figure}

The system (\ref{prey_predator systems}) describes the population dynamics of boyciana-fish system with humans interference which disperse by diffusion in the habitat area $\Omega$ (See Fig.\ref{ecologysys}).
 $x_1(t,\textbf{z}),x_2(t,\textbf{z})$ represent the population densities of fish and boyciana at time $t>0$ and spatial position $\textbf{z}\in\Omega$ respectively. $x_3(t,\textbf{z})$ stands for the human density. The Neumann boundary condition
  $$\frac{\partial x_1}{\partial \textbf{n}}
  =\frac{\partial x_2}{\partial \textbf{n}}
  =\frac{\partial x_3}{\partial \textbf{n}}
  =0, \textbf{z}\in \partial\Omega$$
  means that (\ref{prey_predator systems}) is self-contained and has no population flux across the boundary $\partial\Omega$, so that
$\partial\Omega$ acts as a perfect barrier to dispersal.
 The interaction between fish and boyciana is based on two ratio-dependent functional response functions $$\frac{cx_1x_2}{x_1+\alpha x_2}, \frac{mx_1x_2}{x_1+\alpha x_2}$$
 where $c$ is the capturing rate (or catching efficiency) of the boyciana, $m$ is the conversion rate. Fish population follows the logistic growth in the absence of boyciana and human. $d$ is the death rate of boyciana.

When the distribution of the individuals is not uniform and depends on different spatial locations, the standard method to describe the spatial effects is to introduce the diffusion terms. Therefore the diffusion coefficient matrix about the fish and the boyciana is introduced, that is
$$\left(
         \begin{array}{cc}
           d_{1} & 0 \\
           0& d_{2} \\
         \end{array}
       \right)\Delta x$$
where ${x}=[x_1(t,\textbf{z}),x_2(t,\textbf{z})]^{T}$, $\Delta=\Sigma_{i=1}^{2}\partial^{2}_{z_i}$ is the Laplace operator. It is well known that the appearance of the spatial dispersal makes the dynamics and behaviors of the boyciana-fish system even more complicated.

In particular, for a wetland ecosystem the influence of human can be regarded as an invasive species and not be affected by other species. Thus in (\ref{prey_predator systems}), $h_1$ and $h_2$ are human's interference coefficients on fish $x_1$ and boyciana $x_2$ respectively. Therefore, mathematically the influence of human activities (for example the economic interest) on fish and boyciana are represented by $-h_1x_3$ and $-h_2x_3$  which we add in the first two equations of (\ref{prey_predator systems}).

The third equation of (\ref{prey_predator systems}) is derived from the well-known Fisher equation
\begin{equation}\label{human_distr01}
\frac{\partial x_3(t,\textbf{z})}{\partial t}=\Delta x_3(t,\textbf{z})+rx_3(t,\textbf{z})(1-x_3(t,\textbf{z})),\textbf{z}\in\Omega,t>0.
\end{equation}
where the nonlinear function $rx_3(1-x_3)$ referred to as a logistic nonlinearity. Since the local human population distribution can reach a time independent dynamic balance in a short time. Thus (\ref{human_distr01}) degenerates to  the following elliptic equation
\begin{equation}\label{human_distr02}
0=\Delta x_3+r_3x_3(t,\textbf{z})(1-x_3(t,\textbf{z})),\textbf{z}\in\Omega.
\end{equation}

\begin{rem}

 The aim of our work is to propose some qualitative analysis about three species (boyciana, fish and human) interrelated spatiotemporal ecological wetland model expressed in term of PDAEs. And as application, we illustrate some numerical simulation and prediction to show the effective of our results.
\end{rem}
The remaining part of this paper is organized as follows. In the next section, we analyze the singular part of the system (\ref{prey_predator systems}) including the solution's existence property, bounded property and eigenvalue estimation. In section 4, we analyze the local asymptotic stability of the system's positive equilibrium by using the qualitative theory of dynamical system. And with the prior energy estimation theory we discuss the global stability property of the system (\ref{prey_predator systems}). Finally, in section 5, under some reasonable assumptions, following the qualitative results discussed in Section 2-4, we propose some simulation examples on the system (\ref{prey_predator systems}) to illustrate the effectiveness of our theoretical results. On the other hand, with the data base we collected from the wetland conservation, we carry some numerical analysis and model parameters fitting to predict the boyciana population in the future.

Notations: $\mathbb{N}_0$ is natural number set. $\Omega$ is a bounded plane domain with the boundary $\partial\Omega$. $\|\cdot\|$ denotes the Euclidean norm for vectors. For a symmetric matrix $M$, $M> (<) 0$ means that it is positive (negative) definite. $I$ is the identity matrix. The superscript $T$ is used for the transpose. Matrices, if not explicitly stated, are assumed to have compatible dimensions. For the convenience, we define the following Hilbert space:
  $$\textrm{H}_2(\Omega)\triangleq \{x:\Omega\times[0,+\infty)\rightarrow \mathbb{R}^n \ \textrm{and} \ \|x\|<\infty\}$$
  with inner product and $ L_{2}$-norm  respectively defined by
$$
\langle x,y\rangle \triangleq
\int_{\Omega}x^{T}y\textrm{d}\textbf{z},
\|x\| \triangleq
\{\int_{\Omega}\|x\|^2\textrm{d}\textbf{z}\}^{1/2}.
$$

\section{Stability of the positive steady state}
\subsection{Studied on the human distribution subsystem}
In this section, we investigate the properties of the subsystem in (\ref{prey_predator systems}). That is
\begin{equation}\label{elliptic_subsystem}
\left\{
\begin{aligned}
  &0 =\Delta x_3(\textbf{z})+rx_3(\textbf{z})(1-x_3(\textbf{z})), \textbf{z}\in\Omega,\\
  &\frac{\partial x_3(\textbf{z})}{\partial \textbf{n}}=0,\textbf{z}\in\partial\Omega.
\end{aligned}
\right.
\end{equation}
The solution of (\ref{elliptic_subsystem}) represents the spatial distribution of the humans population which have direct influence on the populations of boyciana and their food. With the existing results in elliptic PDE theory \cite{1972Sat,Smoller1994} we know that the solution of  the system (\ref{elliptic_subsystem}) is not unique. Specially, $x_3(\textbf{z})$ depends not only on the choice of parameter $r$ but also on the shape of the domain $\Omega$ \cite{Smoller1994}. Now, we discuss the existence and uniqueness of the solution  of (\ref{elliptic_subsystem}) and give some ecological description on the positive solution.

First of all, it is obvious that the system (\ref{elliptic_subsystem}) has a coupled upper and lower constant solutions $x^u_3=1,x^l_3=0$ in $\Omega$. The trivial solution $x^u_3=1$ can be interpreted as the humans population density reaches the wetland ecology system capacity limit. In other words, the wetland area is facing the humans  overdevelopment threat. On the other hand, the solution $x^l_3=0$ means that the wetland ecology system is in a human-free environment.

Under the realistic circumstance, the humans population density has decreasing property along $z_1$ axis (see Fig.\ref{ecologysys}) which is perpendicular to the coastline. Therefore what we are interested in is the nontrivial solution of Neumann problem (\ref{elliptic_subsystem}). We refer to \cite{1972Sat} for  discussing the existence of nontrivial positive solution of (\ref{elliptic_subsystem}).

 Let $\underline{x}_3=0$, $\bar{x}_3=1$. Noticing that $\bar{x}_3, \underline{x}_3$ are a couple of upper and lower solutions of (\ref{elliptic_subsystem}) and satisfy the given assumptions in \cite{1972Sat}. By directly using the comparison principle we know that there exists a nonconstant positive solution $x_3^p(\textbf{z})$ satisfies $0\leq x_3^p\leq 1$.

For the uniqueness of the nontrivial positive solutions, \cite[Section 3.5.3]{2003Cantrell} shows that (\ref{reduced_prey_predator systems}) has the unique positive solution $x_3^p$ if and only if
$$r>\lambda_1(\Omega)$$
where $\lambda_1(\Omega)$ is the first eigenvalue of the following elliptic problem
\begin{equation}\label{Eig-Prob}
\left\{
\begin{aligned}
  &0 =\Delta x_3(\textbf{z})+rx_3(\textbf{z}), \textbf{z}\in\Omega,\\
  &\frac{\partial x_3(\textbf{z})}{\partial \textbf{n}}=0,\textbf{z}\in\partial\Omega.
\end{aligned}
\right.
\end{equation}
In addition, the solution $x_3$ will approach zero as $r\searrow\lambda_1(\Omega)$. And by the variational
formula the following inequality holds
$$\int_{\Omega}x_3^p \textrm{d}\textbf{z}\leq (r-\lambda_1(\Omega))|\Omega|,$$
where $|\Omega|$ is the area of $\Omega$. The unique positive solution $x_3^p$ of (\ref{elliptic_subsystem}) has the property that $x_3^p\rightarrow 1$ uniformly on every closed subset of $\Omega$ as $r\rightarrow\infty$.
In summary, we have the following result
\begin{thm}[Uniqueness of positive solution]\label{unique solut thm}
let $\lambda_1(\Omega)$ be the first eigenvalue of (\ref{Eig-Prob}). If $r>\lambda_1(\Omega)$ holds, then (\ref{elliptic_subsystem}) has an unique nontrivial positive solution $x_3^p$ such that $0\leq x_3^p\leq 1$ and
\begin{equation}\label{unique solut thm_eqn01}
  \lim_{r\searrow\lambda_1(\Omega)}x_3^p=0,\ \ \int_{\Omega}x_3^p \textrm{d} \textbf{z}\leq (r-\lambda_1(\Omega))|\Omega|.
\end{equation}
Additionally, if $\partial\Omega\in C^1$ then
\begin{equation}\label{unique solut thm_eqn02}
  \lim_{r\rightarrow\infty} x_3^p=1,\ \ \textbf{z}\in\overline{\Omega'}.
\end{equation}
where $\overline{\Omega'}\subset\Omega$ is closed in $\Omega$.
\end{thm}
\begin{rem}
By the above discussion we obtain three nonnegative solutions of (\ref{elliptic_subsystem}) $\underline{x}_3=0, \bar{x}_3=1$ and $x_3^p$. These correspond to three different kinds of population distribution: human-free, population limit and normal non-uniform distribution. The biological interpretation of Thm.\ref{unique solut thm} is that if the wetland system has the humans carrying capacity $1$ then over all of $\Omega$ the total population would be $|\Omega|$. So the integral inequality (\ref{unique solut thm_eqn01}) shows that the total population reduces with $r$ approach $\lambda_1(\Omega)$. The limit equality (\ref{unique solut thm_eqn02}) shows if the spatial scale of $\Omega$ as measured by $r$ is sufficiently large then the population density in $\Omega$ will be close to its carrying capacity $x_3^u$ on $\Omega$ except for a relatively narrow strip near the boundary $\partial\Omega$. To guarantee the existence and uniqueness of the nontrivial positive human population distribution, we assume $r>\lambda_1(\Omega)$ throughout this study.
\end{rem}
\subsection{Local asymptotic stability}
In this section, we will focus on the steady state solutions of (\ref{prey_predator systems}). Since the semi-linear elliptic subsystem (\ref{elliptic_subsystem}) is independent of time. Thus, we will be mainly concerned with the following parabolic subsystem:
\begin{equation}\label{parabolic system}
\left\{
\begin{aligned}
  &\frac{\partial x_1}{\partial t}
  =d_{1}\Delta x_1+x_1\left(1-x_1-\frac{cx_2}{x_1+\alpha x_2}- h_1x_3\right), \textbf{z}\in\Omega, t>0,\\
  &\frac{\partial x_2}{\partial t}
  =d_{2}\Delta x_2+x_2\left(-d+\frac{mx_1}{x_1+\alpha x_2}-h_2x_3\right), \textbf{z}\in\Omega, t>0,\\
   &\frac{\partial x_1}{\partial \textbf{n}}
  =\frac{\partial x_2}{\partial \textbf{n}}=0, \textbf{z}\in\partial\Omega, t>0,\\
  &x_1(0,\textbf{z})=x_1^0(\textbf{z})\geq 0,x_2(0,\textbf{z})=x_2^0(\textbf{z})\geq 0, \textbf{z}\in\Omega.
\end{aligned}
\right.
\end{equation}
By substituting the humans distribution $x_3=0,1$ and $x_3^p$ into (\ref{parabolic system}), we have
\begin{equation}\label{parabolic system01}
\left\{
\begin{aligned}
  &\frac{\partial x_1}{\partial t}
  =d_{1}\Delta x_1+x_1\left(1-x_1-\frac{cx_2}{x_1+\alpha x_2}\right), \textbf{z}\in\Omega, t>0,\\
  &\frac{\partial x_2}{\partial t}
  =d_{2}\Delta x_2+x_2\left(-d+\frac{mx_1}{x_1+\alpha x_2}\right), \textbf{z}\in\Omega, t>0,\\
   &\frac{\partial x_1}{\partial \textbf{n}}
  =\frac{\partial x_2}{\partial \textbf{n}}=0, \textbf{z}\in\partial\Omega, t>0,\\
  &x_1(0,\textbf{z})=x_1^0(\textbf{z})\geq 0,x_2(0,\textbf{z})=x_2^0(\textbf{z})\geq 0, \textbf{z}\in\Omega.
\end{aligned}
\right.
\end{equation}
\begin{equation}\label{parabolic system02}
\left\{
\begin{aligned}
  &\frac{\partial x_1}{\partial t}
  =d_{1}\Delta x_1+x_1\left(1-x_1-\frac{cx_2}{x_1+\alpha x_2}- h_1\right), \textbf{z}\in\Omega, t>0,\\
  &\frac{\partial x_2}{\partial t}
  =d_{2}\Delta x_2+x_2\left(-d+\frac{mx_1}{x_1+\alpha x_2}-h_2\right), \textbf{z}\in\Omega, t>0,\\
   &\frac{\partial x_1}{\partial \textbf{n}}
  =\frac{\partial x_2}{\partial \textbf{n}}=0, \textbf{z}\in\partial\Omega, t>0,\\
  &x_1(0,\textbf{z})=x_1^0(\textbf{z})\geq 0,x_2(0,\textbf{z})=x_2^0(\textbf{z})\geq 0, \textbf{z}\in\Omega.
\end{aligned}
\right.
\end{equation}
\begin{equation}\label{parabolic system03}
\left\{
\begin{aligned}
  &\frac{\partial x_1}{\partial t}
  =d_{1}\Delta x_1+x_1\left(1-x_1-\frac{cx_2}{x_1+\alpha x_2}- h_1x_3^p(\textbf{z})\right), \textbf{z}\in\Omega, t>0,\\
  &\frac{\partial x_2}{\partial t}
  =d_{2}\Delta x_2+x_2\left(-d+\frac{mx_1}{x_1+\alpha x_2}-h_2x_3^p(\textbf{z})\right), \textbf{z}\in\Omega, t>0,\\
   &\frac{\partial x_1}{\partial \textbf{n}}
  =\frac{\partial x_2}{\partial \textbf{n}}=0, \textbf{z}\in\partial\Omega, t>0,\\
  &x_1(0,\textbf{z})=x_1^0(\textbf{z})\geq 0,x_2(0,\textbf{z})=x_2^0(\textbf{z})\geq 0, \textbf{z}\in\Omega.
\end{aligned}
\right.
\end{equation}
Firstly, from the ecological point of view the above three parabolic type equations systems describe three different kinds of ecological situations. The system (\ref{parabolic system01}) represents that the wetland system is in the original ecological environment situation which is not realistic for Beidaihe wetland system as a international tourism.
For the system (\ref{parabolic system02}), it represents the wetland system is in the maximize development of resources situation. According to the local birds research study, the main threaten factor to the bird diversity in Beidaihe wetland now is from the man-made destruction of the habitation and the over-exploitation of the resource.
In more realistic situation, human, boyciana and the wetland fish are coexistence in wetland environment which is described by the system (\ref{parabolic system03}). Now, mathematically we discuss the stability properties of the systems (\ref{parabolic system01})-(\ref{parabolic system03}) respectively. And study the properties relations among three ecological system.

\subsubsection{Local asymptotic stability of (\ref{parabolic system01})}
It is obvious that the positive constant equilibrium of the system (\ref{parabolic system01}) is the positive solution of the following nonlinear equations
\begin{equation}\label{equil_eqn01}
\left\{
\begin{aligned}
  &f_1(x_1,x_2):=x_1\left(1-x_1-\frac{cx_2}{x_1+\alpha x_2}\right)=0,\\
  &f_2(x_1,x_2):=x_2\left(-d+\frac{mx_1} {x_1+\alpha x_2}\right)=0.
\end{aligned}
\right.
\end{equation}
By some simple analysis for (\ref{equil_eqn01}), we have
\begin{thm}
 When $1-\frac{\alpha}{c}<\frac{d}{m}<1$, the system (\ref{parabolic system01}) has a unique positive equilibrium
$e_1:=(u_1^\ast, \upsilon_1^\ast)$ with
$$ u_1^\ast=1-\frac{c(m-d)}{m\alpha}=1-\frac{c}{\alpha}(1-\frac{d}{m}), \upsilon_1^\ast=\frac{m(1-u_1^\ast)u_1^\ast}{cd}=\frac{1}{\alpha}(\frac{m}{d}-1)u_1^\ast$$
In addition, $0<u_1^\ast<1, \upsilon_1^\ast> 0$; $u_1^\ast\geq 1$ and $\upsilon_1^\ast\leq 0$ when $\frac{d}{m}\geq 1$.
\end{thm}
Now, with the linearized method of dynamic system and the eigenvalue theory of PDE we discuss the local asymptotic stability of the positive constant steady state $e_1$.
Considering the linearizing system of (\ref{parabolic system01}) at $e_1$
\begin{equation}\label{parabolic system01 Linearized}
  \frac{\partial x}{\partial t}=\mathscr{L}x
\end{equation}
where the linear operator $\mathscr{L}$ be defined by
\begin{equation}\label{JACB01}
  \mathscr{L}=D_1\Delta+\bar{f}_J
\end{equation}
with $\bar{f}_J$ is the Jacobian matrix of $\bar{f}=(f_1,f_2)^T$ at $e_1$. Let
$$0=\lambda_0<\lambda_1<\cdots<\lambda_n<\cdots$$
be the eigenvalues of the operator $-\Delta$ on $\Omega$ with Neumann BC.
The corresponding eigenfunctions are represented by $\phi_n (n\in \mathbb{N}_0)$.
Thus, $\lambda_n, \phi_n (n\in \mathbb{N}_0)$ satisfy
\begin{equation}\label{SL problem}
\left\{
\begin{aligned}
      &-\Delta\phi=\lambda\phi \ \ \textrm{in}\ \Omega,\\
      &\frac{\partial\phi}{\partial \textbf{n}}=0 \ \ \textrm{on} \ \partial\Omega. \\
\end{aligned}
\right.
\end{equation}
Then the function sequence $\{\phi_n\}_{n=0}^{\infty}$ forms an orthonormal base of $\textrm{H}^2(\Omega)$. It should be noticed that the eigenvector w.r.t. $\lambda_0=0$ is $\phi_0=\textrm{const}$. The corresponding solution is trivial which can not influence the stability of the system. Therefore, we are concerned in the following infinity dimensional ODE systems (see \cite{Jiang2015} for detail):
\begin{equation}\label{parabolic system01 Linearized infty}
 \dot{X}_n=\mathscr{L}_nX_n, n\in\mathbb{N}_0.
\end{equation}
where $X_n=\langle x, \phi_n\rangle=\int_\Omega x(t,\textbf{z})\phi(\textbf{z})\textrm{d}\textbf{z},$
\begin{equation}\label{JACB02}
  \mathscr{L}_n=-D_1\lambda_n+\bar{f}_J
  \end{equation}
substituted $\bar{f}_J$ with $e_1$, one get
\begin{equation}\label{JACB03}
  \mathscr{L}_n=
  \left(
    \begin{array}{cc}
      -\frac{m^2\alpha+cd^2-cm^2+d_1\lambda_n\alpha m^2}{m^2\alpha} & -c\frac{d^2}{m^2} \\
      \frac{(m-d)^2}{m\alpha} & -\frac{\lambda_nd_2m+md-d^2}{m}\\
    \end{array}
  \right)
\end{equation}
Let $\textrm{Dt}(e_1), \textrm{Tr}(e_1)$ denote the determinant and the trace of the matrix $\mathscr{L}_n$ respectively. Then
\begin{equation}\label{Determ and Trace01}
\left\{
\begin{aligned}
  \textrm{Dt}(e_1)=&\frac{m^2d\alpha-cd^3-md^2\alpha-m^2dc+2md^2c}{m^2\alpha}\\
+&\frac{\lambda_nd^2cd_2+\lambda_nm^2\alpha d_2-\lambda_nm^2cd_2+\lambda_n^2m^2\alpha d_1d_2+\lambda m^2\alpha d_1d-\lambda m d^2 \alpha d_1}{m^2\alpha},\\
  \textrm{Tr}(e_1)=&
-\frac{m^2\alpha\lambda_n(d_1+d_2)+m^2\alpha+d^2c-m^2c}{m^2\alpha}.
\end{aligned}
\right.
\end{equation}
Noticing that $m>d$, by the monotonicity of eigenvalues $\{\lambda_n\}_{n=1}^{\infty}$, through some directly computing we have the following conclusion
\begin{thm}\label{equil stable thm01}
Under the condition $1-\frac{\alpha}{c}<\frac{d}{m}<1$,\\
\textrm{S1}) When $c\leq\alpha$, the positive equilibrium $e_1$ is locally asymptotically stable for any diffusion coefficients $d_1, d_2>0.$\\
\textrm{S2}) When $c>\alpha$, the positive equilibrium $e_1$ is locally asymptotically stable, if and only if
$$\lambda_1(d_1+d_2)-\frac{c}{\alpha}\frac{m-d}{m}\frac{d}{m}>0$$
where $\lambda_1$ is the first eigenvalue of S-L problem (\ref{SL problem}).
\end{thm}
\begin{rem}
Define the domains $\textrm{I}, \textrm{II}$ and $\textrm{III}$ in the $(c,\alpha, m, d)$ space respectively by
\begin{eqnarray}
  & & \textrm{I}=\{c\leq \alpha, 0<\frac{d}{m}<1\}, \\
  & & \textrm{II}=\{c> \alpha, 1-\frac{\alpha}{c}<\frac{d}{m}<1, \frac{c}{\alpha}\frac{m-d}{m}\frac{d}{m}\leq \lambda_1(d_1+d_2)\},\\
  & & \textrm{III}=\{c> \alpha, 1-\frac{\alpha}{c}<\frac{d}{m}<1, \frac{c}{\alpha}\frac{m-d}{m}\frac{d}{m}> \lambda_1(d_1+d_2)\}.
\end{eqnarray}
we assert that if $(c,\alpha, m, d)\in \textrm{I,II}$ then the positive equilibrium $e_1$ of system (\ref{parabolic system01}) is locally asymptotically stable, else if $(c,\alpha, m, d)$ is in $\textrm{III}$ then $e_1$ is unstable (See subsection \ref{Sim_section_human_free} for the numerical illustration).

From a biological point of view, Thm. \ref{equil stable thm01} implies that if there is no human interference ($x_3=0$), when the capturing rate $c$ is low with the high conversion rate $m$, the Ciconia boyciana population $x_2$ and their fish food $x_1$ can maintain the stable coexistence situation. In addition, the diffusion term $\lambda_1(d_1+d_2)$ in the inequality shows that better diffusion effect helps maintain the balance of the ecosystem. The above described circumstances consistent with experience.
\end{rem}
\subsubsection{Local asymptotic stability of (\ref{parabolic system02})}
Analogously, for the system (\ref{parabolic system02}), the positive constant equilibrium of  is the positive solution of the following nonlinear equations
\begin{equation}\label{equil_eqn02}
\left\{
\begin{aligned}
  &x_1\left(1-x_1-\frac{cx_2}{x_1+\alpha x_2}-h_1\right)=0,\\
  &x_2\left(-d+\frac{mx_1} {x_1+\alpha x_2}-h_2\right)=0.
\end{aligned}
\right.
\end{equation}
By some simple analysis for (\ref{equil_eqn02}), we have
\begin{thm}
 When $1-\frac{\alpha}{c}(1-h_1)<\frac{d+h_2}{m}<1$, the system (\ref{parabolic system02}) has a unique positive equilibrium
$e_2:=(u_2^\ast, \upsilon_2^\ast)$ with
$$ u_2^\ast=1-h_1-\frac{c(m-d-h_2)}{m\alpha}=1-h_1-\frac{c}{\alpha}(1-\frac{d+h_2}{m}),$$
$$\upsilon_2^\ast=\frac{m(1-h_1-u_2^\ast)u_2^\ast}{c(d+h_2)}=\frac{1}{\alpha}(\frac{m}{d+h_2}-1)u_2^\ast.$$
In addition, $0<u_2^\ast<1-h_1, \upsilon_2^\ast> 0$; $u_2^\ast\geq 1$ and $\upsilon_2^\ast\leq 0$ when $\frac{d+h_2}{m}\geq 1$.
\end{thm}
In a similar manner, the local asymptotically stability at $e_2$ can be investigated by linearization method.
 The concerned infinity dimensional ODEs system matrix is
\begin{equation}\label{JACB04}
  \mathscr{L}^\ast_n=
  \left(
    \begin{array}{cc}
      -\frac{m^2\alpha+cd^2-cm^2+d_1\lambda_n\alpha m^2}{m^2\alpha}-h_1 & -c\frac{d^2}{m^2} \\
      \frac{(m-d)^2}{m\alpha} & -\frac{\lambda_nd_2m+md-d^2}{m}-h_2\\
    \end{array}
  \right)
\end{equation}
Noticing that $\mathscr{L}^\ast_n=\mathscr{L}_n-\lambda_n\textrm{diag}(h_1, h_2)$,
if we denote the determinant and the trace of the matrix $\mathscr{L}^\ast_n$ as $\textrm{Dt}(u_2^\ast), \textrm{Tr}(u_2^\ast)$ respectively, then
\begin{equation}\label{Determ and Trace02}
\left\{
\begin{aligned}
  \textrm{Dt}(u_2^\ast)=&\textrm{Dt}(u_1^\ast)-(a_{11}\lambda_nh_1+ a_{22}\lambda_nh_2)+h_1h_2\lambda_n^2,\\
  \textrm{Tr}(u_2^\ast)=& \textrm{Tr}(u_2^\ast)-(h_1+h_2)\lambda_n.
\end{aligned}
\right.
\end{equation}
where $a_{11},a_{22}$ are the diagonal elements of $\textrm{Dt}(u_1^\ast)$. From Thm.\ref{equil stable thm01} we know that
$$a_{11}<0, a_{22}<0,$$
if $\lambda_1(d_1+d_2)-\frac{c}{\alpha}\frac{m-d}{m}\frac{d}{m}>0 $ is provided.

Thus by the nonnegativity of $\lambda_n$, one can get that $$\textrm{Dt}(u_2^\ast)>0, \textrm{Tr}(u_2^\ast)<0.$$
By summarizing, we have
\begin{thm}\label{equil stable thm02}
Under the condition $1-\frac{\alpha}{c}(1-h_1)<\frac{d+h_2}{m}<1$,\\
\textrm{S3}) When $c\leq\alpha$, the positive equilibrium $e_2$ is locally asymptotically stable for any diffusion coefficients $d_1, d_2>0$ and interference coefficients $h_1, h_2>0$.\\
\textrm{S4}) When $c>\alpha$, the positive equilibrium $e_2$ is locally asymptotically stable, if
$$\lambda_1(d_1+d_2)-\frac{c}{\alpha}\frac{m-d}{m}\frac{d}{m}>0$$
is provided. Here $\lambda_1$ is the first eigenvalue of S-L problem (\ref{SL problem}).
\end{thm}
\begin{rem}
From Thm. \ref{equil stable thm02}, we get the diffusion coefficients $d_1,d_2$ also determine the stability of the equilibrium $e_2$. Additionally, the interference coefficients $h_1, h_2$ makes the equilibrium $e_2$ more smaller, even change to zeros. If we assume $h_2>\frac{m\alpha}{c}h_1$, then the positive condition $1-\frac{\alpha}{c}(1-h_1)<\frac{d+h_2}{m}<1$ implies $1-\frac{\alpha}{c}<\frac{d}{m}<1$, not vice versa. It is means that the steady state $e_1$ can become unstable with human influence (See subsection \ref{Sim_section_overdevelop} for the numerical illustration).
\end{rem}
\subsection{Global Stability of the nontrivial steady state of (\ref{parabolic system02})}
In the last two subsections, the linearized method can solve the local stability of the constant steady state $e_1, e_2$. However, for the nontrivial steady state $e_3$, it is the function of the spatial variable $\textbf{z}$. Obviously, under this situation we can't expect a similar situation in front of the equilibrium point.

We first obtain the following global attractor result for the solution of (\ref{parabolic system03}) which can be similarly shown as in \cite{2013Ko} by the simple comparison argument for parabolic equations.
\begin{thm}
If $1-\frac{\alpha}{c}<\frac{d}{m}<1$, $1-\frac{\alpha}{c}(1-h_1)<\frac{d+h_2}{m}<1$, then for all the nonnegative solution of (\ref{parabolic system03}) $x(t, \textbf{z}),\textbf{z}\in\overline{\Omega}$, the following inequalities stand
\begin{eqnarray}\label{Eqn_Global_attractor01}
\label{Eqn_Global_attractor0101}
  & &\limsup_{t\rightarrow\infty} x_1(t,\textbf{z})\leq \bar{x}_1:=1,\\
  \label{Eqn_Global_attractor0102}
  & &\limsup_{t\rightarrow\infty} x_2(t,\textbf{z})\leq \bar{x}_2:=\frac{1}{\alpha}\frac{m-d}{d},\\
  \label{Eqn_Global_attractor0103}
  & &\liminf_{t\rightarrow\infty} x_1(t,\textbf{z})\geq \underline{x}_1:=1-h_1-\frac{c}{\alpha},\\
  \label{Eqn_Global_attractor0104}
  & &\liminf_{t\rightarrow\infty} x_2(t,\textbf{z})\geq \underline{x}_2:=
  \frac{1}{\alpha}\frac{m-(d+h_2)}{d+h_2}(1-h_1-\frac{c}{\alpha}).
\end{eqnarray}
Consequently, the domain given by
$\mathscr{A}:=[1-h_1-\frac{c}{\alpha}, 1]\times
[\frac{1}{\alpha}\frac{m-(d+h_2)}{d+h_2}(1-h_1-\frac{c}{\alpha}), \frac{1}{\alpha}(\frac{m-d}{d})]$
is a positively invariant region for global solutions of system (\ref{parabolic system03}).
\end{thm}
\textbf{Proof:} Let $x_1, x_2$ be a solution of (\ref{parabolic system03}). Then from the first equation of (\ref{parabolic system03}) one can observe that $x_1$ satisfies
 \begin{equation}\label{Eqn_Global_attractor02}
\left\{
\begin{aligned}
  &\frac{\partial x_1}{\partial t}
  \leq d_{1}\Delta x_1+x_1(1-x_1), \textbf{z}\in\Omega, t>0,\\
  &\frac{\partial x_1}{\partial \textbf{n}}=0, \textbf{z}\in\partial\Omega, t>0,\\
  &x_1(0,\textbf{z})=x_1^0(\textbf{z})\geq 0, \textbf{z}\in\Omega.
\end{aligned}
\right.
\end{equation}
In view of the comparison principle of parabolic equation, one can get that
$$\limsup_{t\rightarrow\infty} x_1(t,\textbf{z})\leq 1.$$
Thus the second equation of (\ref{parabolic system03}) yields that $x_2$ satisfies
\begin{equation}\label{Eqn_Global_attractor03}
\left\{
\begin{aligned}
  &\frac{\partial x_2}{\partial t}
  \leq d_{2}\Delta x_2+x_2(-d+\frac{m}{1+\alpha x_2}), \textbf{z}\in\Omega, t>0,\\
   &\frac{\partial x_2}{\partial \textbf{n}}=0, \textbf{z}\in\partial\Omega, t>0,\\
  &x_2(0,\textbf{z})=x_2^0(\textbf{z})\geq 0, \textbf{z}\in\Omega.
\end{aligned}
\right.
\end{equation}
 Noticing that $\bar{x}_2=\frac{1}{\alpha}\frac{m-d}{d}$ satisfies
 \begin{equation}\label{Eqn_Global_attractor04}
\left\{
\begin{aligned}
  &\frac{\partial \bar{x}_2}{\partial t}
  = d_{2}\Delta \bar{x}_2+\bar{x}_2(-d+\frac{m}{1+\alpha \bar{x}_2}), \textbf{z}\in\Omega, t>0,\\
   &\frac{\partial \bar{x}_2}{\partial \textbf{n}}=0, \textbf{z}\in\partial\Omega, t>0.
\end{aligned}
\right.
\end{equation}
Hence the inequality (\ref{Eqn_Global_attractor0102}) holds.

In addition, from the first equation of (\ref{parabolic system03}) we can also get
 \begin{equation}\label{Eqn_Global_attractor05}
\left\{
\begin{aligned}
  &\frac{\partial x_1}{\partial t}
  \geq d_{1}\Delta x_1+x_1(1-\frac{c}{\alpha}-h_1-x_1), \textbf{z}\in\Omega, t>0,\\
  &\frac{\partial x_1}{\partial \textbf{n}}=0, \textbf{z}\in\partial\Omega, t>0,\\
  &x_1(0,\textbf{z})=x_1^0(\textbf{z})\geq 0, \textbf{z}\in\Omega.
\end{aligned}
\right.
\end{equation}
It follows that the inequality (\ref{Eqn_Global_attractor0103}) holds. From the second equation of (\ref{parabolic system03}) we have
\begin{equation}\label{Eqn_Global_attractor06}
\left\{
\begin{aligned}
  &\frac{\partial x_2}{\partial t}
  \geq d_{2}\Delta x_2+x_2(-d+\frac{m\underline{x}_1}{\underline{x}_1+\alpha x_2}-h_2), \textbf{z}\in\Omega, t>0,\\
   &\frac{\partial x_2}{\partial \textbf{n}}=0, \textbf{z}\in\partial\Omega, t>0,\\
  &x_2(0,\textbf{z})=x_2^0(\textbf{z})\geq 0, \textbf{z}\in\Omega.
\end{aligned}
\right.
\end{equation}
Therefore, (\ref{Eqn_Global_attractor0104}) holds.
\begin{rem}
For the system (\ref{parabolic system03}), the reaction function vector is
\begin{eqnarray}\label{reaction term}
\bar{f}(x)=\left(
                    \begin{array}{c}
                      x_1 \left( 1-x_1-\frac{cx_2}{x_1+\alpha x_2}- h_1x_3^p(\textbf{z})\right) \\
                      x_2\left(-d+\frac{mx_1}{x_1+\alpha x_2}-h_2x_3^p(\textbf{z})\right)
                    \end{array}
                  \right).
\end{eqnarray}
Since the Jacobian matrix of $\bar{f}$ is
\begin{equation}\label{reaction_term_Jacobi}
  L=\left(
    \begin{array}{cc}
      1-c\alpha\left(\frac{x_2}{x_1+\alpha x_2}\right)^2-2x_1-h_1x_3^p & -c\left(\frac{x_1}{x_1+\alpha x_2}\right)^2 \\
      m\alpha\left(\frac{x_2}{x_+\alpha x_2}\right)^2 & m\left(\frac{x_1}{x_1+\alpha x_2}\right)^2-d-h_2x_3^p \\
    \end{array}
  \right),
\end{equation}
 $\bar{f}$ is a mixed quasi-monotone function vector in $\mathbb{R}^2_+=\{(x_1, x_2)|x_1\geq 0, x_2\geq 0\}$. Therefore the above theorem implies that $\bar{x}:=(\bar{x}_1,\bar{x}_2)$ and $\underline{x}:=(\underline{x}_1,\underline{x}_2)$ are a pair of coupled upper and lower solutions of (\ref{parabolic system03}). Consequently, by\cite[Charpter 8 Thm. 3.3]{2012Pao} there exists a solution $x(t,\textbf{z})$ of (\ref{parabolic system03}) with $\underline{x}\leq x(t,\textbf{z})\leq \bar{x}$.
  \end{rem}
\begin{rem}
From the first two inequalities (\ref{Eqn_Global_attractor0101}),(\ref{Eqn_Global_attractor0102})
of the above theorem, we observe that if $m<d$, then $\lim_{t\rightarrow\infty}x_2=0$ uniformly on $\overline{\Omega}$. Ecologically, the boyciana population will tends toward extinction.
The last two inequalities (\ref{Eqn_Global_attractor0103}),(\ref{Eqn_Global_attractor0104}) give sufficient conditions such that the positive solution of (\ref{parabolic system02}) has the persistence property. That is, we provide some necessary conditions on parameters such that the boyciana and fish  always coexist with humans influence:
$$h_1<1-\frac{c}{\alpha}, h_2<m-d.$$
This shows that it is reasonable to expect the persistence of boyciana and fish when there is a suitable weak humans influence.
\end{rem}
In the following, we propose some exponential stability property on the system (\ref{prey_predator systems}) by PDAEs energy estimation theory. For the following considerations, it will be simplest to assume homogeneous Neumann BC though our theoretical result is applicative to the Dirchlet BC case.
\begin{lem}[Poincare inequality \cite{Smoller1994}] \label{lem_Poincare}
Let $x\in W^{2}_{2}(\Omega)$, then if $\mu_1$ is the smallest positive eigenvalue of $-\Delta$ on $\Omega$ (with the appropriate boundary conditions) the following Poincar\'{e} inequalities hold:
\begin{equation}\label{Poincare inequalities01}
  \|\nabla x\|^2\geq \mu_1 \|x-\bar{x}\|^2, \|\Delta x\|^2\geq \mu_1\|\nabla x\|^{2}\ if\ \frac{dx}{d\textbf{n}}=0\ on \ \partial\Omega,
\end{equation}
\begin{equation}\label{Poincare inequalities02}
  \|\nabla x\|^2\geq \mu_1 \|x\|^2 \ if\ x=0\ on \ \partial\Omega,
\end{equation}
where $\bar{x}=\frac{1}{|\Omega|}\int_{\Omega}xd\textbf{z}.$
\end{lem}
\begin{thm}\label{Thm Energy estim}
Assume that $\breve{x}(t,\textbf{z})=(x_1(t,\textbf{z}),x_2(t,\textbf{z}))^T$ is a bounded solution of (\ref{parabolic system03}). Assume that spectral radius of the Jacobian matrix (\ref{JACB03}) is $\rho$, $\mu_1$ is the smallest positive eigenvalue of $-\Delta$ on $\Omega$, $d_m=\min\{d_1, d_2\}$ and $$\delta=d_m\mu_1-(\rho+h_M)>0.$$
where $h_M=\max\{h_1, h_2\}$, then
\begin{equation}\label{Thm Energy estim eqn01}
\int_{\Omega}\sum^{2}_{i=1}\frac{\partial {\breve{x}}^{T}}{\partial z_i}\frac{\partial {\breve{x}}}{\partial z_i}\textrm{d}\textbf{z}=\int_{\Omega}\|\nabla \breve{x}\|^2d\textbf{z} \leq c_1\exp(-\delta t)+c_2,
\end{equation}
\begin{equation}\label{Thm Energy estim eqn02}
\int_{\Omega}\|{\breve{x}}(t,\textbf{z})-\breve{x}_M(t)\|\textrm{d}\textbf{z}\leq c_3\exp(-\delta t)+c_4
\end{equation}
hold for positive constant $c_1, c_2$ and $c_3, c_4$. Here, $$\breve{x}_M(t)=\frac{1}{|\Omega|}\int_{\Omega}\breve{x}(t,\textbf{z})\textrm{d}\textbf{z}$$
is the spatial average function.
\end{thm}
\textbf{Proof:} Let us introduce the energy integral (Lyapunov function)
\begin{equation}\label{Thm Energy estim Energy function}
  E(t)=\frac{1}{2}\int_{\Omega}\|\nabla \breve{x}\|^2\textrm{d}\textbf{z}
  =\frac{1}{2}\int_{\Omega}\sum^{2}_{i=1}
  \left[
  \left(\frac{\partial {x_i}}{\partial z_1}\right)^2
  +\left(\frac{\partial {x_i}}{\partial z_2}\right)^2
  \right]
  \textrm{d}\textbf{z}
\end{equation}
By computing the derivative of $E(t)$, one get
\begin{eqnarray}
  \frac{\textrm{d}}{\textrm{d}t}E(t) &=&\int_{\Omega}\sum^{2}_{i=1}
  \left[
  \frac{\partial {x_i}}{\partial z_1}\frac{\partial^2 {x_i}}{\partial z_1\partial t}
  +\frac{\partial {x_i}}{\partial z_2}\frac{\partial^2 {x_i}}{\partial z_2\partial t}
    \right]
    \textrm{d}\textbf{z}
\end{eqnarray}
Noticing that $\breve{x}$ satisfies the system (\ref{prey_predator systems}), it follows from above
\begin{eqnarray}
  \frac{\textrm{d}}{\textrm{d}t}E(t) &=&\int_{\Omega}\sum^{2}_{i=1}
  \left(
      \frac{\partial {\breve{x}}^{T}}{\partial z_i}\frac{\partial}{\partial z_i}( D_1\Delta {x}+\bar{f}(x))
      \right)
      \textrm{d}\textbf{z}
\end{eqnarray}
where $D_1=\textrm{diag}(d_1, d_2)$, $\bar{f}=(f_1, f_2)^T$ is the nonlinear reaction vector function (\ref{reaction term}). Applying for divergence theorem, we get
 \begin{eqnarray}
  \frac{\textrm{d}}{\textrm{d}t}E(t) =  -\int_{\Omega}\Delta{\breve{x}}^{T}D_1\Delta{\breve{x}}\textrm{d}\textbf{z}
   +\int_{\Omega}\sum^{2}_{i=1}
      \frac{\partial {\breve{x}}^{T}}{\partial z_i}\frac{\partial \bar{f}}{\partial z_i}\textrm{d}\textbf{z}=I_1+I_2.
\end{eqnarray}
For the first integral $I_1$,  the following estimation holds
\begin{equation}\label{Thm Energy estim eqn03}
I_1\leq -d_{m}\int_{\Omega}\|\Delta\breve{x}\|^2\textrm{d}\textbf{z}
\leq-d_{m}\mu_1\int_{\Omega}\|\nabla\breve{x}\|^2\textrm{d}\textbf{z}.
\end{equation}
The second integral $I_2$ can be estimated as
\begin{equation}\label{Thm Energy estim eqn04}
  I_2=
  \int_{\Omega}\sum^{2}_{j=1}\sum^{2}_{i=1}\left\{
  \frac{\partial {x_j}}{\partial z_i}
  \left(
  \frac{\partial {f_j}}{\partial x_1}\frac{\partial {x_1}}{\partial z_i}
  +\frac{\partial {f_j}}{\partial x_2}\frac{\partial {x_2}}{\partial z_i}
  +\frac{\partial {f_j}}{\partial x_3}\frac{\partial {x_3}}{\partial z_i}
  \right)
  \right\}\textrm{d}\textbf{z}
\end{equation}
or
\begin{equation}\label{Thm Energy estim eqn05}
  I_2=
  \int_{\Omega}\nabla\breve{x}^{T}
  \frac{\partial {\bar{f}}}{\partial \breve{x}}
  \nabla\breve{x}\textrm{d}\textbf{z}
  +\int_{\Omega}\left(
  \frac{\partial {f_1}}{\partial x_3}\nabla x_1\cdot\nabla x_3
  +\frac{\partial {f_2}}{\partial x_3}\nabla x_2\cdot\nabla x_3
  \right)
  \textrm{d}\textbf{z}
\end{equation}
 where $\frac{\partial {\bar{f}}}{\partial \breve{x}}$ is the Jacobi matrix $\bar{f}$ w.r.t. $\breve{x}$. Thus
\begin{equation}\label{Thm Energy estim eqn06}
  I_2\leq
  \rho\int_{\Omega}\|\nabla\breve{x}\|^2\textrm{d}\textbf{z}
  +h_1\int_{\Omega}\|\nabla x_1\|^2\textrm{d}\textbf{z}
  +h_2\int_{\Omega}\|\nabla x_2\|^2\textrm{d}\textbf{z}
  +(h_1+h_2)\int_{\Omega}\|\nabla x_3\|^2\textrm{d}\textbf{z}.
\end{equation}
 According to Lemma \ref{lem_Poincare}, (\ref{Thm Energy estim eqn03})  and (\ref{Thm Energy estim eqn06}) imply
\begin{equation}\label{Thm Energy estim eqn07}
  \frac{\textrm{d}}{\textrm{d}t} E(t)\leq
  (-d_m \mu_1+\rho+h_M) E(t)+ (h_1+h_2)\int_{\Omega}\|\nabla x_3\|^2\textrm{d}\textbf{z}.
\end{equation}
Noticing that
$$\Delta x_3(\textbf{z})+rx_3(\textbf{z})(1-x_3(\textbf{z}))=0, \textbf{z}\in\Omega$$
stands with the Neumann boundary condition, we have
\begin{equation}\label{Thm Energy estim eqn08}
  \frac{\textrm{d}}{\textrm{d}t} E(t)\leq
  (-d_m \mu_1+\rho+h_M) E(t)+ r(h_1+h_2)\int_{\Omega} x_3^2(1-x_3)\textrm{d}\textbf{z}.
\end{equation}
From the assumption $ d_m\mu_1>\rho+h_M$, with Gronwall inquality we obtain that
\begin{equation}\label{Thm Energy estim eqn09}
  E(t)\leq (E(0)-q)\exp( -\delta t)+q,
\end{equation}
where $\delta=d_m\mu_1-(\rho+h_M)> 0$,
$$q=\frac{r(h_1+h_2)\int_{\Omega} x_3^2(1-x_3)\textrm{d}\textbf{z}}{\delta}.$$
Therefore, (\ref{Thm Energy estim eqn01}) holds. By lemma \ref{lem_Poincare}, (\ref{Thm Energy estim eqn01}) implies (\ref{Thm Energy estim eqn02}).
\begin{rem}
Our proposed result generalized the energy estimation result on common parabolic type PDEs system. In \cite{2013Ko} the energy function $E(t)$ tends to zero as $t\rightarrow\infty$. This assertion is not right in our system (\ref{prey_predator systems}). However, from above theorem, the value of $E(t)$ is asymptotically decreasing tends to a small value $q$ when
 $$\frac{r(h_1+h_2)\int_{\Omega} x_3^2(1-x_3)\textrm{d}\textbf{z}}{d_m\mu_1-(\rho+h_M)}\rightarrow 0.$$
 The corresponding conditions to the above inequality are
 $$d_m\rightarrow\infty\ \textrm{and}\ \rho, h_M, r, h_1, h_2\rightarrow 0.$$
 If we consider two special cases $x_3=0$ and $x_3=1$, then $q=0$ and the system (\ref{parabolic system03}) translate into
 (\ref{parabolic system01}) and (\ref{parabolic system01}). Immediately, we get the following corollary.
 \end{rem}
\begin{cor}\label{Corl Energy estim01}
Under the conditions described in Thm. \ref{Thm Energy estim}, If $\breve{x}$, $\breve{y}$ are bounded solutions of (\ref{parabolic system01}) and (\ref{parabolic system02}) respectively,
then the following conclusions hold
\begin{equation}\label{Corl Energy estim eqn01}
\int_{\Omega}\sum^{2}_{i=1}\frac{\partial {\breve{x}}^{T}}{\partial z_i}\frac{\partial {\breve{x}}}{\partial z_i}\textrm{d}\textbf{z}=\int_{\Omega}\|\nabla \breve{x}\|^2d\textbf{z} \leq c_1\exp(-\delta_1 t),
\end{equation}
\begin{equation}\label{Corl Energy estim eqn02}
\int_{\Omega}\|{\breve{x}}-\breve{x}_M(t)\|\textrm{d}\textbf{z}\leq c_2\exp(-\delta_1 t)
\end{equation}
\begin{equation}\label{Corl Energy estim eqn03}
\int_{\Omega}\sum^{2}_{i=1}\frac{\partial {\breve{y}}^{T}}{\partial z_i}\frac{\partial {\breve{y}}}{\partial z_i}\textrm{d}\textbf{z}=\int_{\Omega}\|\nabla \breve{y}\|^2d\textbf{z} \leq c_3\exp(-\delta_2 t),
\end{equation}
\begin{equation}\label{Corl Energy estim eqn04}
\int_{\Omega}\|{\breve{y}}-\breve{y}_M(t)\|\textrm{d}\textbf{z}\leq c_4\exp(-\delta_2 t)
\end{equation}
where $c_1, c_2, c_3, c_4$ are all independent positive constants; $\breve{x}_M(t)$, $\breve{y}_M(t)$ are the spatial average functions of $\breve{x}(t)$ and $\breve{y}(t)$ respectively.
 i.e. the state variable vector $\breve{x}$, $\breve{y}$ generated by the system (\ref{parabolic system01}) and (\ref{parabolic system02}) are all exponentially stable and asymptotically converge to their spatial average respectively.
\end{cor}
\begin{rem}
Ecologically, $E(t)$ represents the global spatial mobility integral of boyciana and fish species. It's value can be seen as the combination of the oscillation amplitude about two species populations in the domain. In case of human interference, Thm.\ref{Thm Energy estim} shows that in order to avoid the inference of human behaviors (for example: human fishing, water pollution), boyciana and wetland fish must carry on the unceasing migration. This means the wetland system is in an unstable state. However, if we reduced the human distribution density (set up human restrict area) and enhance the spatial diffusion capacity (for example,avoid river pollution), the wetland ecological system can tend to a stable state (See subsection 5.3 for the numerical illustration).
\end{rem}
\section{Model fitting on real boyciana-fish data}
In this section, as application of our boyciana-fish-human model, we evaluate the performance of the proposed PDAEs model by comparing the density calculated by the model with the actual bird data set. We first deal with the observation data. The original observation data will be rearranged as one dimensional location data. Then the one dimension space PDAEs model is built. And we use the least square fitting technique provided in \cite{2013Lei_WangHY} to search for parameters to best fit the real data.
\subsection{Real data treatment}\label{data_treatment}
With the help of qinhuangdao BRBS (Bird Reserve and Banding Station), this real data is collected  from six different locations at different times (from 2001 to 2014). Tab.\ref{Table_boyciana_distrib} shows the discrete density distribution of boyciana $\tilde{x}_2(t,z_i) (i=1,\cdots,6)$ at six bird-watching locations in normalization format.
\begin{center}
\begin{table}
  \centering
  \caption{Boyciana density distribution $\tilde{x}_2(t,z)$ in six observatories.}\label{Table_boyciana_distrib}
  \begin{tabular}{c|c|c|c|c}
  \hline
  Time (year)     & 2001 & 2002 & $\cdots$ & 2014  \\ \hline
  $z_1$ (Loc.1) & 0.001136 & 0.001047 & $\cdots$ & 0.000717  \\
  $z_2$ (Loc.2) & 0.001324 & 0.001150 & $\cdots$ & 0.000981  \\
  $z_3$ (Loc.3) & 0.000919 & 0.000802 & $\cdots$ & 0.000742  \\
  $z_4$ (Loc.4) & 0.000946 & 0.000909 & $\cdots$ & 0.000766  \\
  $z_5$ (Loc.5) & 0.001100 & 0.000966 & $\cdots$ & 0.000396  \\
  $z_6$ (Loc.6) & 0.000900 & 0.000866 & $\cdots$ & 0.000266  \\
  \hline
\end{tabular}
\end{table}
\end{center}
One can find that the density value at Location 1 is higher than the other five locations. It is worth notices that this observatory is located in the forest where few people enter and is far away from the coastline. This confirms that boyciana depend on and prefer a region without interference or less interference for building their nest. Now, since these locations have different linear distances from the Beidaihe coastline, accordingly we discrete the spatial axis $z$ as six points $z_1,\cdots, z_6$
which correspond to the six bird observatories locations.

For the fish density distribution $\tilde{x}_1(t,z)$, it should be noticed that the "fish density" represents the fish food available quantity. Considering the range of boyciana foraging activities is at about $0.5$ to $5$ km, we assume that the distance from the boyciana's habitat to the coastline represents degree of different fish food catching. That is the quantity of available fish food linearly negative dependent on the distance to the coastline. Thus, in Tab.\ref{Table_fish_distrib}, the density distribution at $z_6$ position $\tilde{x}_1(t,z_6)$ is normalized created by the real quantity of fish production of Qinhuangdao in the past 14 years (from 2001 to 2014). And the other positions $\tilde{x}_1(t,z_i) (i=1,\cdots,5)$ is linearly decreasing with respect to the coastline distance $z$.

\begin{center}
\begin{table}
\centering
  \caption{Normalization fish density distribution $\tilde{x}_1(t,z)$.}\label{Table_fish_distrib}
  \begin{tabular}{c|c|c|c|c}
\hline
Time (year)&	2001&	2002&$\cdots$&2014\\ \hline
$z_1$	&0.1&	0.09&$\cdots$& 0.03 \\
$z_2$	&0.12&	0.12& $\cdots$&0.07\\
$z_3$	&0.14&	0.14&$\cdots$&0.07\\
$z_4$	&0.16&	0.15&$\cdots$&0.10\\
$z_5$	&0.18&	0.18&$\cdots$&0.12\\
$z_6$	&0.2&	0.19&$\cdots$&0.14\\ \hline
\end{tabular}
\end{table}
\end{center}
\subsection{Spatial assumption}
We make some reasonable assumptions to reduce the spatial domain of the original model (\ref{prey_predator systems}). From subsection\ref{data_treatment}, we know the populations densities of boyciana and fish are assumed as single space variables. Mathematically, the state variables vector satisfies
$$x(t,z_1,z_2)=x(t,z_1).$$
Then the spatial domain $\Omega$ is longitudinal compression into a 'pipe' line
 $$\Omega=(0, l).$$
 Moreover, the wetland boundary $\partial\Omega$ is changed into two points: $\textbf{z}=0, \textbf{z}=l$ which are both closed without flux under Neumann BC.
For the convenience of writing, we define $x(t,z)=x(t,z_1)$ and choose the spatial domain of system (\ref{prey_predator systems}) as $\Omega=(0,\pi)$ to simplify calculation.

Under above assumptions, (\ref{prey_predator systems}) becomes an 1D-spatial PDAEs as following:
\begin{equation}\label{reduced_prey_predator systems}
\left\{
\begin{aligned}
  &\frac{\partial x_1}{\partial t}
  =d_{1}\frac{\partial^2 x_1}{\partial z^2}+x_1\left(1-x_1-\frac{cx_2}{x_1+\alpha x_2}- h_1x_3\right), z\in(0,\pi),t>0,\\
  &\frac{\partial x_2}{\partial t}
  =d_{2}\frac{\partial^2 x_2}{\partial z^2}+x_2\left(-d+\frac{mx_1}{x_1+\alpha x_2}-h_2x_3\right), z\in(0,\pi),t>0,\\
  &0
  =\frac{\partial^2 x_3}{\partial z^2}+rx_3(1-x_3), z \in (0,\pi),\\
  &\frac{\partial x_1}{\partial z}
  =\frac{\partial x_2}{\partial z}
  =\frac{\partial x_3}{\partial z}=0, z\in\{0,\pi\},t>0,\\
     &x_1(0,z)=x_1^0(z), x_2(0,z)=x_2^0(z), z\in[0,\pi].\\
\end{aligned}
\right.
\end{equation}
In addition, the corresponding S-L problem (\ref{SL problem}) has the eigenvalues and the eigenfunctions:
\begin{equation}\label{Sim_Eigenvalue_function}
\begin{aligned}
\mu_n=n^2, \phi_n(z) =\sqrt{\frac{2}{\pi}}\cos(nz),\ \ n\in \mathbb{N}_0.
\end{aligned}
\end{equation}
Therefore, the first eigenvalue is $\mu_1([0,\pi])=1$.

\subsection{Model fitting with real data}
In this subsection, we will determine the system's parameters with a numerical fitting method and predict the future development trend of boyciana population. First, we estimate the parameter $r$ of the human distribution subsystem in (\ref{reduced_prey_predator systems}) to fit the real human population distribution.
\begin{equation}\label{Sim_elliptic_subsystem}
\left\{
\begin{aligned}
  &0 =\frac{\partial^2 x_3}{\partial z^2}+rx_3(z)(1-x_3(z)), z\in (0, \pi),\\
  &\frac{\partial x_3}{\partial z}|_{z=0}=\frac{\partial x_3}{\partial z}|_{z=\pi}=0.
\end{aligned}
\right.
\end{equation}
Since the nontrivial positive analytical solution of (\ref{Sim_elliptic_subsystem}) is in some implicit form with an integral term, we provide an approximate function to numerical computing. That is,
\begin{equation}\label{human_distr_funtion}
x_3(z)\approx1-\frac{3}{4}\textrm{sech}^2(\frac{\sqrt{r}z}{2}), z\in (0, \pi).
\end{equation}
It should be noticed that the above approximate function holds the properties in Thm.\ref{unique solut thm}.
\begin{figure}
  \centering
  \includegraphics[width=0.6\textwidth]{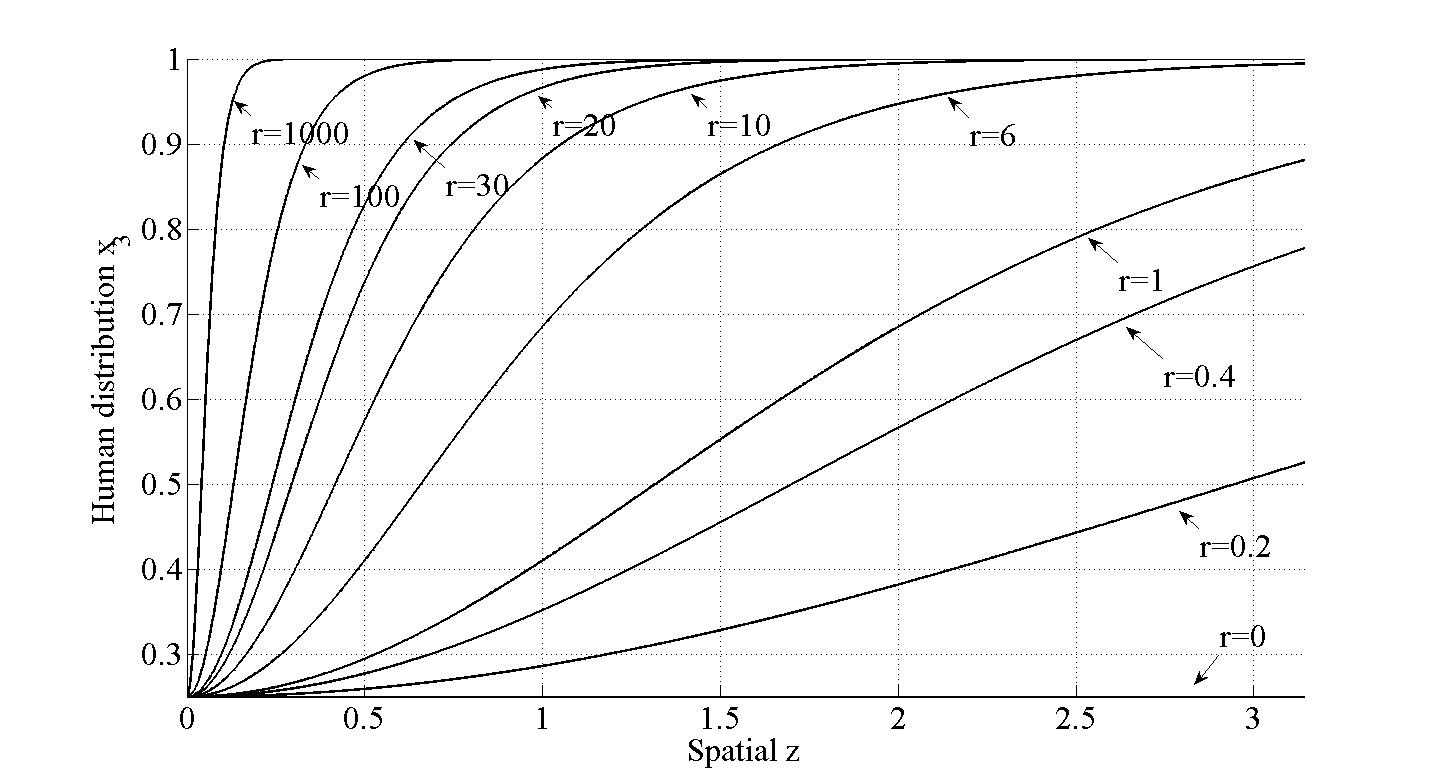}\\
  \caption{Human distribution curve with different parameter values $r$.}\label{Fig_human_distrib}
\end{figure}

\begin{figure}
  \centering
  \includegraphics[width=0.6\textwidth]{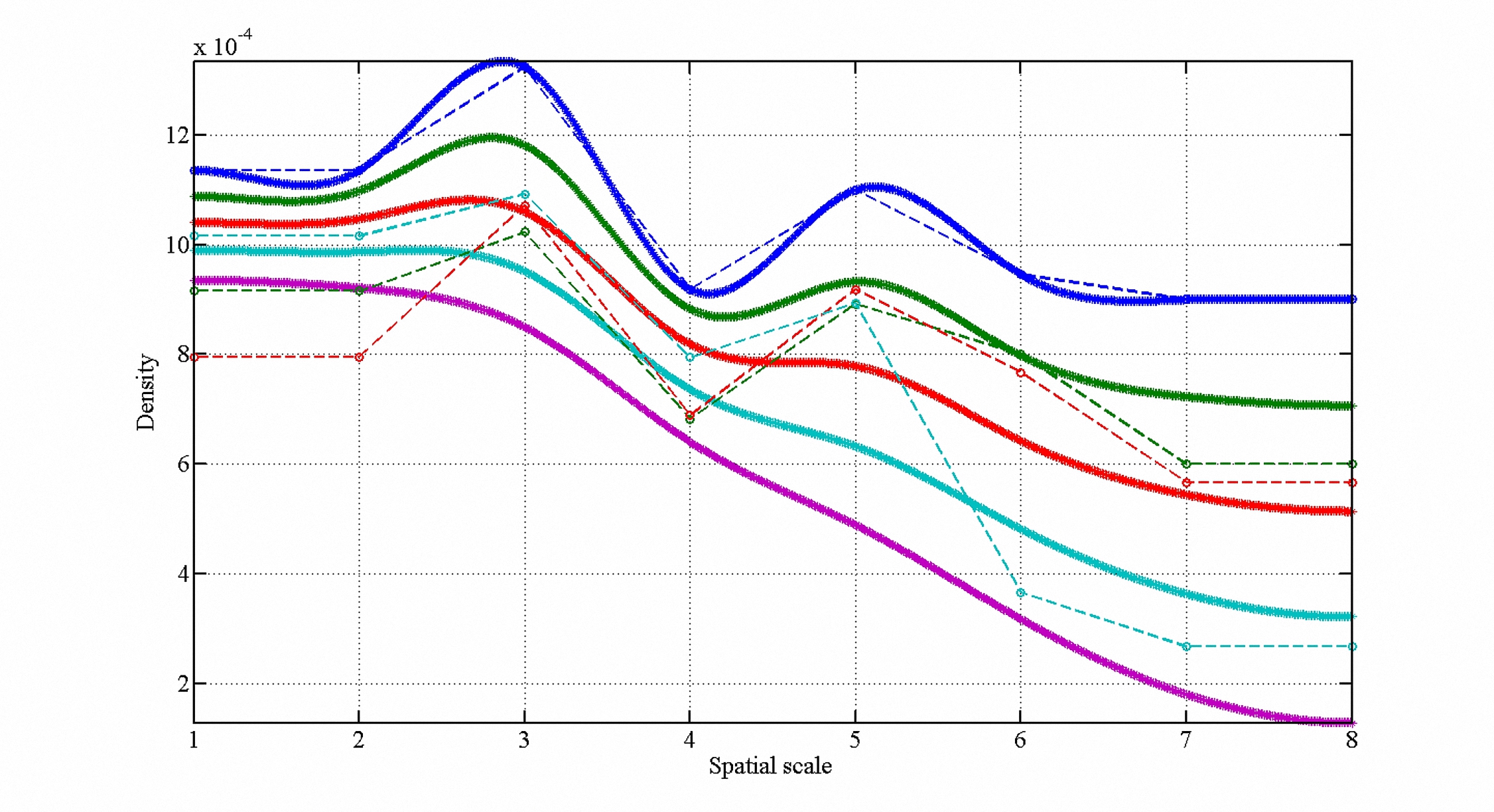}\\
  \caption{Predicted vs. Actual data of boyciana density. The x-axis represents the spatial scale. The y-axis is density. From top to bottom, the first four predicted curves are in years of 2001, 2005, 2009, 2014 respectively.}\label{Fig_Predict_x2}
\end{figure}

Now, we optimize the system's parameter $d_1, d_2, c, m, \alpha, d, h_1, h_2, r$. to fit the practice data. Our optimal problem is
\begin{equation}\label{Optimal_Problem}
\left\{
\begin{aligned}
  &\min\sum_{i=1}^{2}\sum_{j=1}^{14}\sum_{k=1}^{6}\frac{|x_i(t_j,z_k)-\tilde{x}_i(t_j,z_k)|}{|\tilde{x}_i(t_j,z_k)|}\\
  &s.t.\\
  &\frac{\partial x_1}{\partial t}
  =d_{1}\Delta x_1+x_1\left(1-x_1-\frac{cx_2}{x_1+\alpha x_2}- h_1(1-\frac{3}{4}\textrm{sech}^2(\frac{\sqrt{r}z}{2}))\right),\\
  &\frac{\partial x_2}{\partial t}
  =d_{2}\Delta x_2+x_2\left(-d+\frac{mx_1}{x_1+\alpha x_2}-h_2(1-\frac{3}{4}\textrm{sech}^2(\frac{\sqrt{r}z}{2}))\right), \\
    &\frac{\partial x_1}{\partial z}=\frac{\partial x_2}{\partial z}=0,\\
     &x_1(0,z)=\tilde{x}_1(0,z),x_2(0,z)=\tilde{x}_2(0, z),\\
     &d_1, d_2, c, m, \alpha, d, h_1, h_2, r>0.
     \end{aligned}
\right.
\end{equation}
We take the actual data from 2001 as the initial density value of the above optimal problem (\ref{Optimal_Problem}). For the Neumann BC, we use the technique in numerical analysis, the spatial direction we add two sets of data to each side of the observation data. The added data values are the same as their adjacent values. The numerical optimal program is designed with Matlab software.
The initial condition input is
$$d_1=d_2=0.001, d=0.3, c=1, m=1, \alpha=0.5, h_1=0.1, h_2=0.3, r=1.$$
Combining with the boundary value in Tab.\ref{Table_boyciana_distrib} and Tab.\ref{Table_fish_distrib}, we get the numerical optimization result:
$$d_1=0.1185, d_2=0.5773, d=0.1343, c=0.8996, $$
$$m=0.5093, \alpha=0.5535, h_1=1.8102, h_2=0.0172, r=6.0440.$$
Fig.\ref{Fig_Predict_x2}
illustrates the predicted results for boyciana. The dashed lines denote the actual observations for the density in different year, while the starred lines illustrate the density predicted by our PDAEs Model (\ref{Optimal_Problem}). Here, the average prediction accuracy is $95.17\%$. It is effective from the statistical perspective.
\section{Illustrative examples for system stability}
In this section, as an effective demonstration of our previous theoretic analysis, we simulate the steady state for the boyciana-fish-human system
\begin{equation}\label{sim_systems}
\left\{
\begin{aligned}
  &\frac{\partial x_1}{\partial t}
  =d_{1}\frac{\partial^2 x_1}{\partial z^2}+x_1\left(1-x_1-\frac{cx_2}{x_1+\alpha x_2}- h_1x_3\right), z\in(0,\pi),t>0,\\
  &\frac{\partial x_2}{\partial t}
  =d_{2}\frac{\partial^2 x_2}{\partial z^2}+x_2\left(-d+\frac{mx_1}{x_1+\alpha x_2}-h_2x_3\right), z\in(0,\pi),t>0,\\
  &0
  =\frac{\partial^2 x_3}{\partial z^2}+rx_3(1-x_3), z \in (0,\pi),\\
  &\frac{\partial x_1}{\partial z}
  =\frac{\partial x_2}{\partial z}
  =\frac{\partial x_3}{\partial z}=0, z\in\{0,\pi\},t>0.\\
\end{aligned}
\right.
\end{equation}
\subsection{Human free ecological system}\label{Sim_section_human_free}

\begin{figure}
  \centering
  \subfigure[Spatial-temporal profiles]{
    \label{Fig_Sim_Stab_x1_subfig_a} 
    \includegraphics[width=0.5\textwidth]{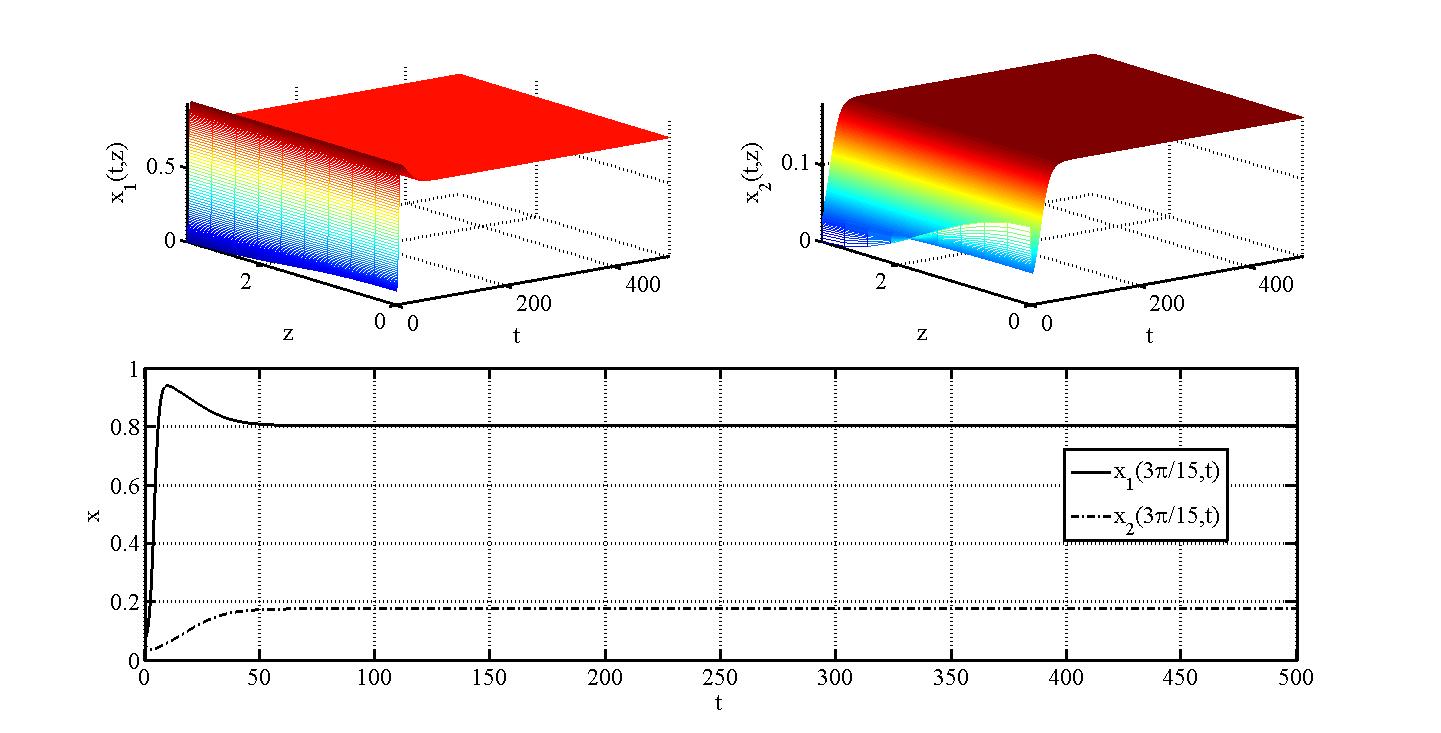}}
  \hspace{-0.05\textwidth}
  \subfigure[Spatial response curve]{
    \label{Fig_Sim_Stab_x1_subfig_b} 
    \includegraphics[width=0.5\textwidth]{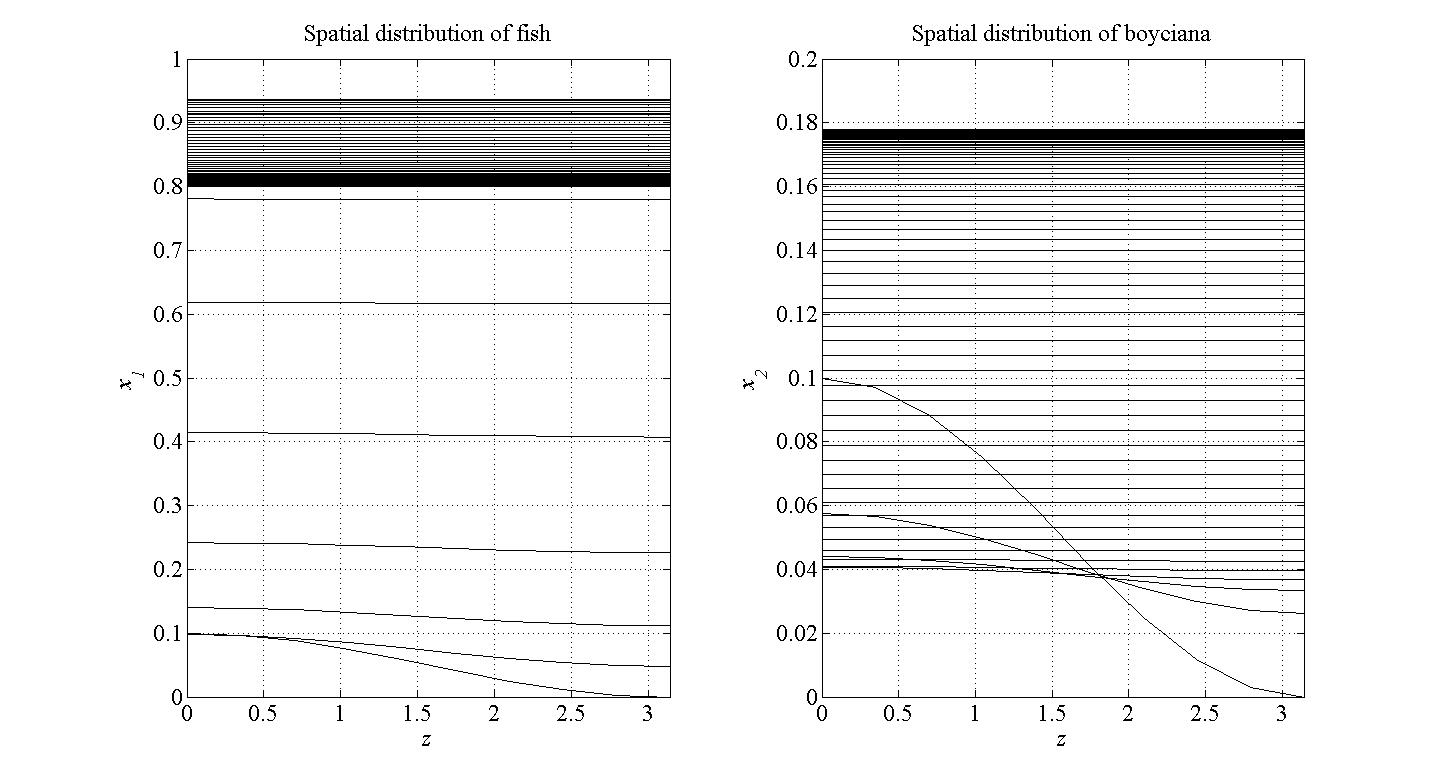}}
  \caption{Numerical results for the local stable property of $e_1$ with $d_1=d_2=1$. (a) On the top part are the spatial-temporal profiles of $x_1(z,t)$ and $x_2(z,t)$. The corresponding time evolution curves is at the bottom part where $z=\frac{\pi}{5}$.
  (b) Spatial response curve about $x_1(z)$ and $x_2(z)$ at differential discrete time points.}
  \label{Fig_Sim_Stab_x1} 
\end{figure}

For the numerical evaluation of the human free system (\ref{sim_systems}), we take the system parameters as
\begin{eqnarray}\label{Sim_Parameter01}
  &c=1.000, \alpha=0.5000, d=0.9000, m=1.000, h_1=0, h_2=0.
\end{eqnarray}
which fulfill the positive condition in Thm.\ref{equil stable thm01}
$$1-\frac{\alpha}{c}<\frac{d}{m}<1.$$
The equilibrium is $e_1\approx(0.8000, 0.1778).$
Simulation results with $d_1=d_2=1.000$ are exemplarily depicted in Fig.\ref{Fig_Sim_Stab_x1}.

\begin{figure}
  \centering
  \subfigure[Spatial-temporal profiles]{
    \label{Fig_Sim_Ustab_x1_subfig_a} 
    \includegraphics[width=0.5\textwidth]{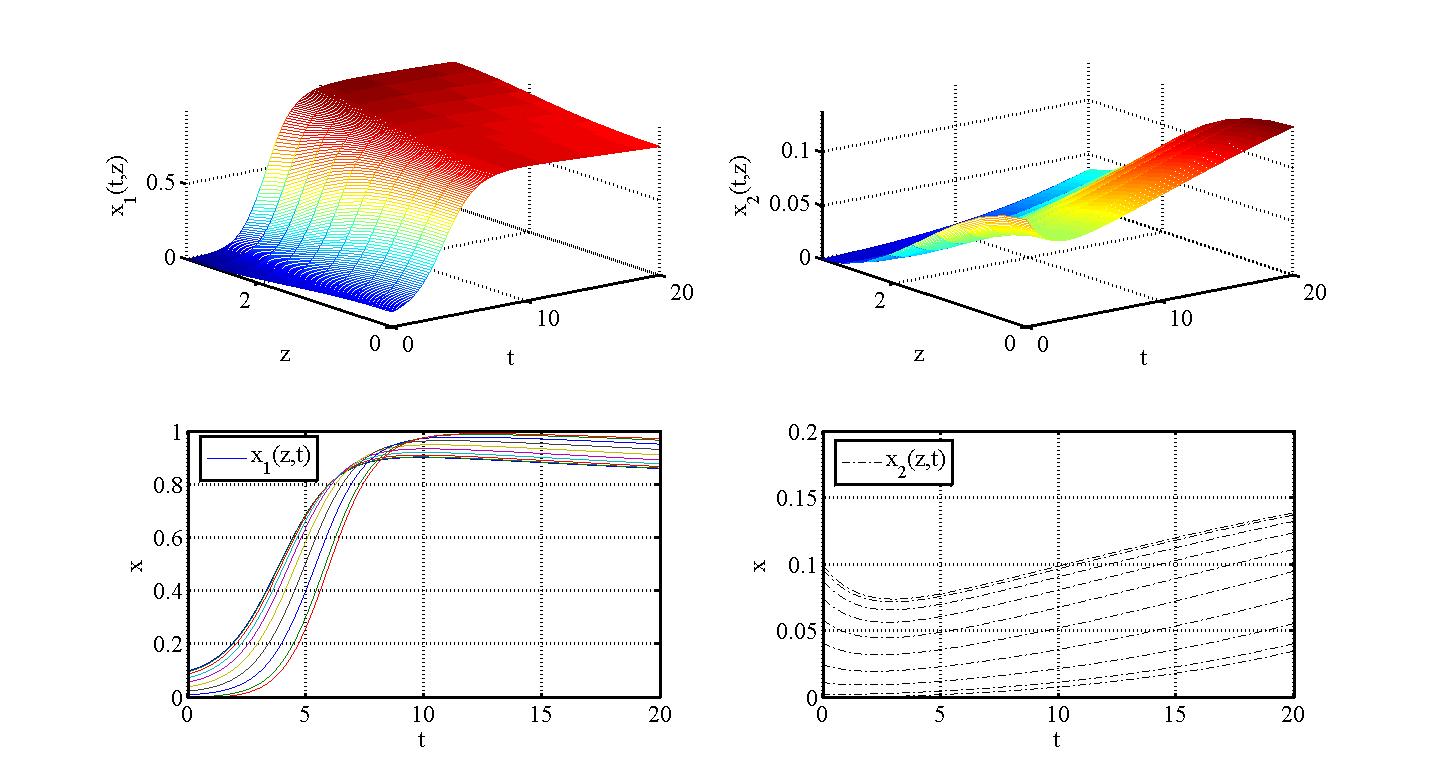}}
  \hspace{-0.05\textwidth}
  \subfigure[Spatial response curve]{
    \label{Fig_Sim_Ustab_x1_subfig_b} 
    \includegraphics[width=0.5\textwidth]{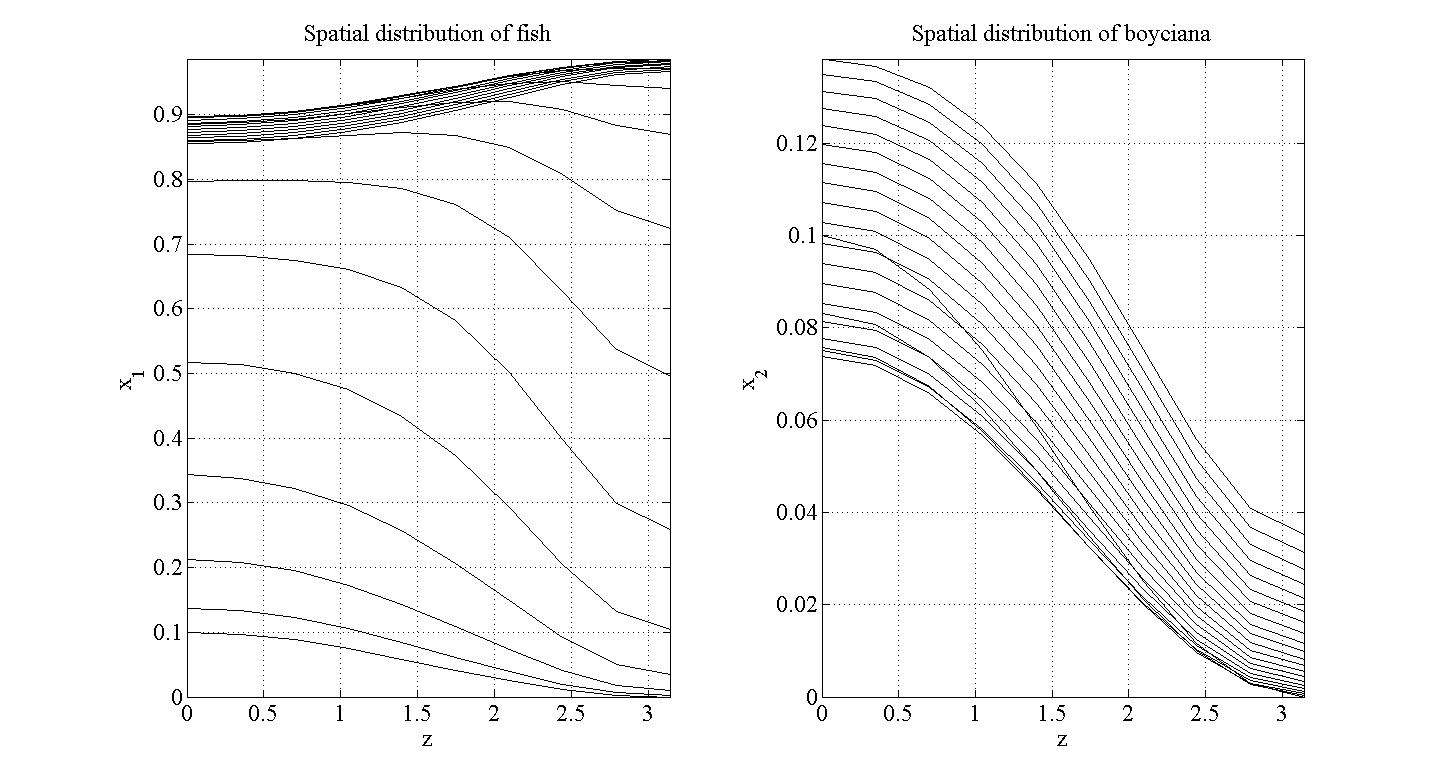}}
  \caption{Numerical results for the local unstable property of $e_1$ with $d_1=d_2=0.01$. (a) On the top part are the spatial-temporal profiles of $x_1(z,t)$ and $x_2(z,t)$. The corresponding time evolution curves is at the bottom part where $z=\frac{\pi}{5}$.
  (b) Spatial response curve about $x_1(z)$ and $x_2(z)$ at differential discrete time points.}
  \label{Fig_Sim_Ustab_x1} 
\end{figure}

In Fig.\ref{Fig_Sim_Stab_x1_subfig_a} the spatiotemporal response surface of fish $x_1(z,t)$ and boyciana $x_2(z,t)$ are depicted on the top part. And also the time evolution of both species at position $z=\frac{\pi}{5}$ is studied on the bottom part.
Fig.\ref{Fig_Sim_Stab_x1_subfig_b} illustrates the spatial distribution curves of boyciana and fish $x_1$, $x_2$ at discrete time points. With these, the equilibrium $e_1$ shows the stable property.On the other hand, in Fig.\ref{Fig_Sim_Ustab_x1}, by decreasing the diffusion coefficients $d_1, d_2$ to $d_1=d_2=0.01$, the numerical results are shown with unstable property.

 Obviously, due to Thm.\ref{equil stable thm01}, if $d_1=d_2=1$ then $$\lambda_1(d_1+d_2)-\frac{c}{\alpha}\frac{m-d}{m}\frac{d}{m} \approx 1.820>0.$$
 $e_1$ is locally stable. If $d_1=d_2=0.01$ then $$\lambda_1(d_1+d_2)-\frac{c}{\alpha}\frac{m-d}{m}\frac{d}{m} \approx -0.160<0.$$
 $e_1$ is unstable. The simulation results correspond with the theoretical conclusion.

\subsection{Overdevelopment ecological system}\label{Sim_section_overdevelop}
\begin{figure}
  \centering
  \subfigure[Spatial-temporal profiles]{
    \label{Fig_Sim_Stab_x2_subfig_a} 
    \includegraphics[width=0.5\textwidth]{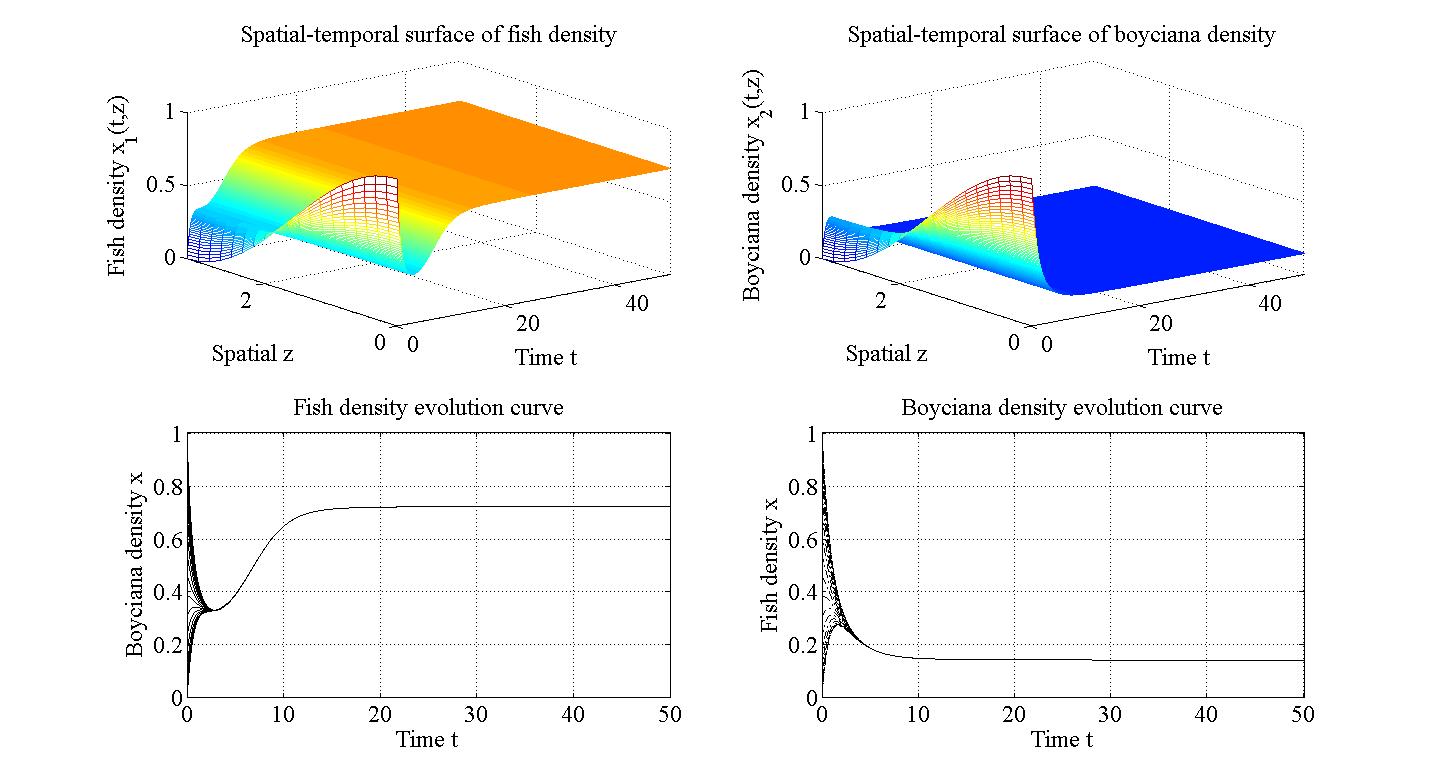}}
  \hspace{-0.05\textwidth}
  \subfigure[Spatial response curve]{
    \label{Fig_Sim_Stab_x2_subfig_b} 
    \includegraphics[width=0.5\textwidth]{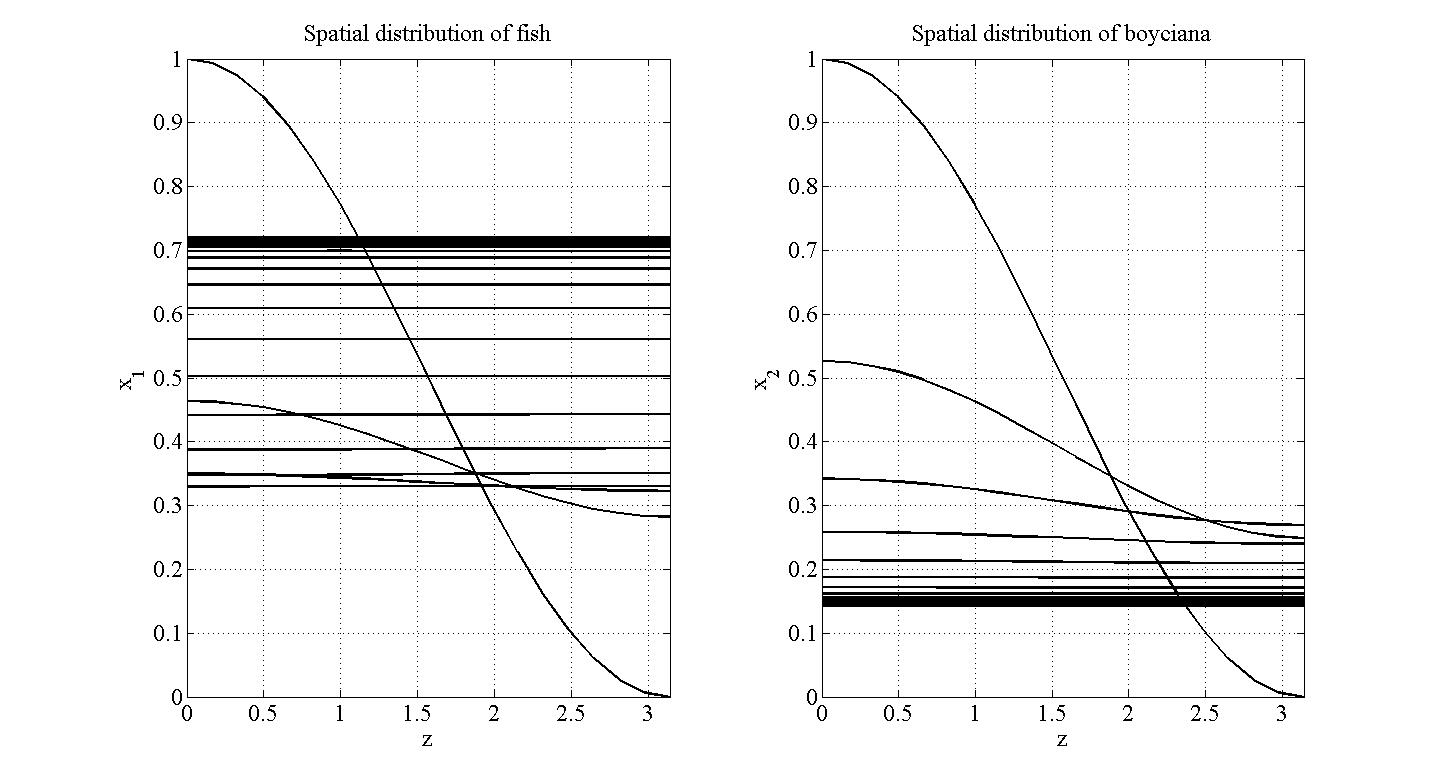}}
  \caption{Numerical results for the local stable property of $e_2$ with $d_1=d_2=1$. (a) Stable steady state profiles and the corresponding time evolution curves at $e_2$. (b) Spatial response curve about $x_1(z)$ and $x_2(z)$ at differential discrete time points.}
  \label{Fig_Sim_Stab_x2} 
\end{figure}

\begin{figure}
  \centering
  \subfigure[Spatial-temporal profiles]{
    \label{Fig_Sim_Ustab_x2_subfig_a} 
    \includegraphics[width=0.5\textwidth]{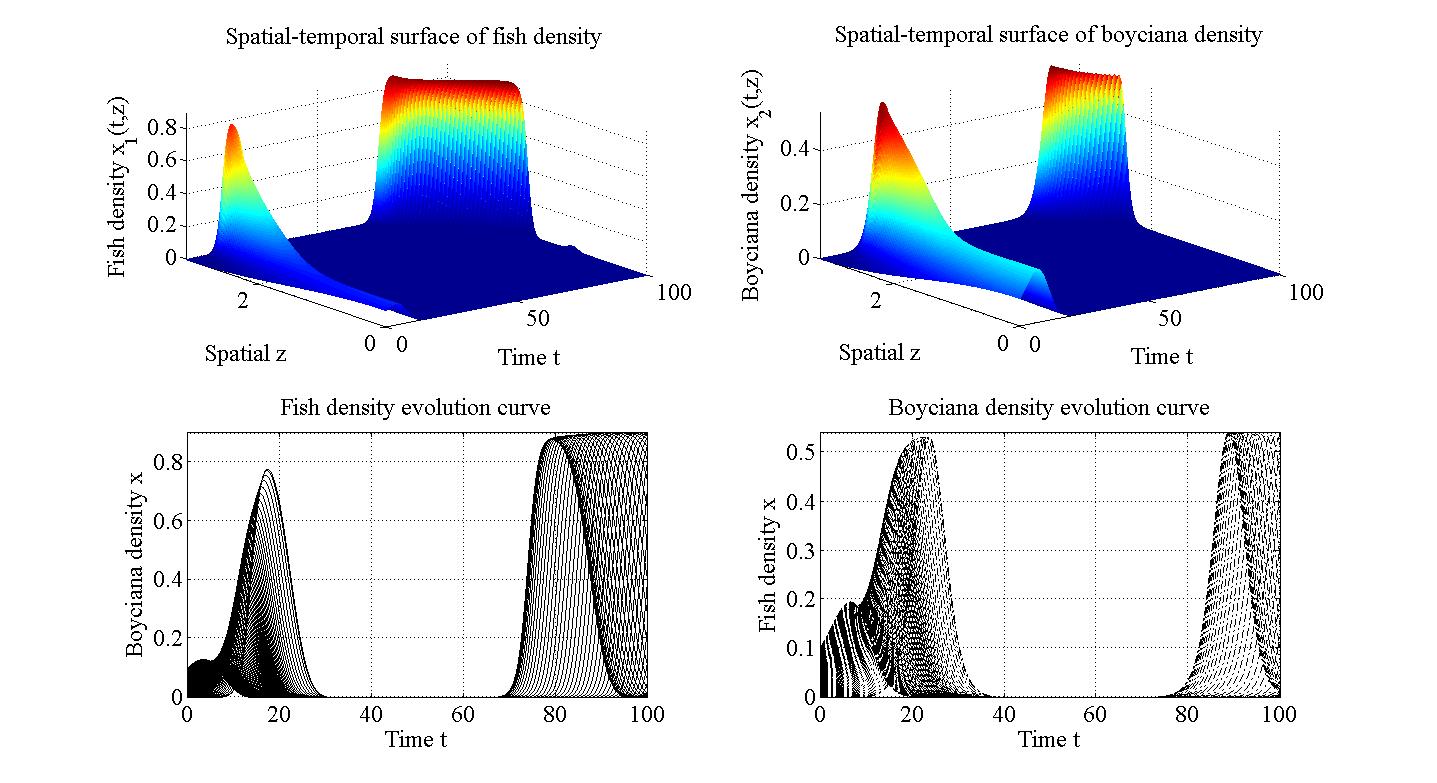}}
  \hspace{-0.05\textwidth}
  \subfigure[Spatial response curve]{
    \label{Fig_Sim_Ustab_x2_subfig_b} 
    \includegraphics[width=0.5\textwidth]{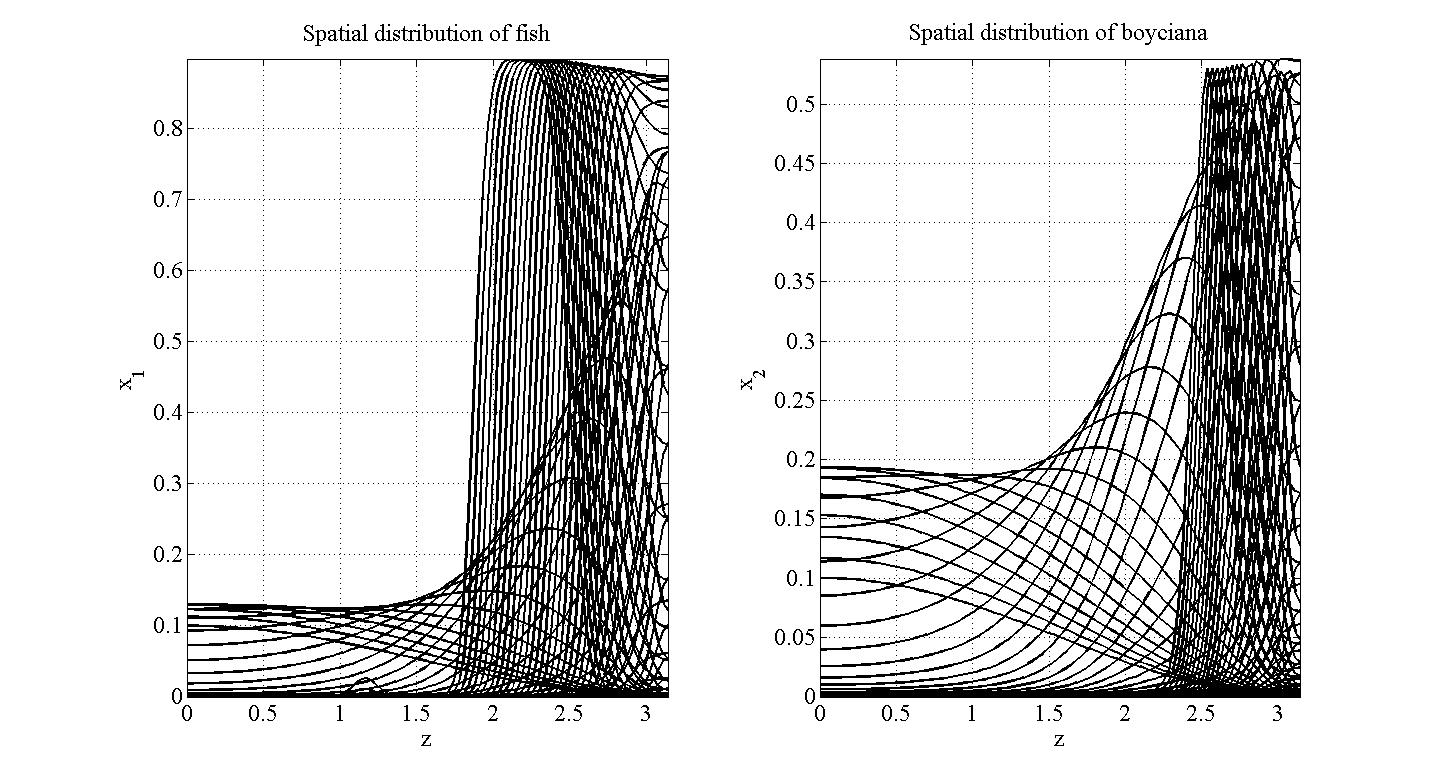}}
  \caption{Numerical results for the local unstable property of $e_2$ with $d_1=d_2=0.001$. (a) Unstable steady state profiles and the corresponding time evolution curves at $e_2$. (b) Spatial response curve about $x_1(z)$ and $x_2(z)$ at differential discrete time points.}
  \label{Fig_Sim_Ustab_x2} 
\end{figure}

For the system (\ref{sim_systems}), considering Thm \ref{equil stable thm02}, following from (\ref{Sim_Parameter01}), we choose the system's parameters as
\begin{eqnarray}\label{Sim_Parameter02}
  &c=1.000, \alpha=0.5000, d=0.9000, m=1.000, h_1=0.1, h_2=0.01.
\end{eqnarray}
where $h_1=0.1$ and $h_2=0.01$ are humans interference coefficients. The other parameters of (\ref{Sim_Parameter02}) are the same as (\ref{Sim_Parameter01}).
By directly computing, one can see that the above parameters fulfill the positive condition
$$1-\frac{\alpha}{c}(1-h_1)<\frac{d+h_2}{m}<1.$$
The equilibrium is $e_2\approx(0.7200, 0.1424).$
Simulation results in the $(z,t)$ domain are exemplarily depicted in Fig.\ref{Fig_Sim_Stab_x2} and Fig.\ref{Fig_Sim_Ustab_x2}.
Fig.\ref{Fig_Sim_Stab_x2} illustrates the stable property of $e_2$. Fig.\ref{Fig_Sim_Ustab_x2} is the unstable case. It is worth noting that in Fig.\ref{Fig_Sim_Stab_x2} the densities of boyciana and fish are both decreasing more than the previous human-free model in Fig.\ref{Fig_Sim_Stab_x1}. This corresponds with the theoretical result.


\subsection{Simulation on boyciana-fish-human coexistence system}
In this case, the nontrivial positive human distribution is introduced in the boyciana-fish model. Therefore we illustrate that even in the case of other conditions unchanged if the humans population density increases the boyciana-fish-human ecological system can turn from stable state into a instability state.

\begin{figure}
  \centering
  \subfigure[Spatial-temporal profiles]{
    \label{Fig_Sim_Stab_x3_subfig_a} 
    \includegraphics[width=0.5\textwidth]{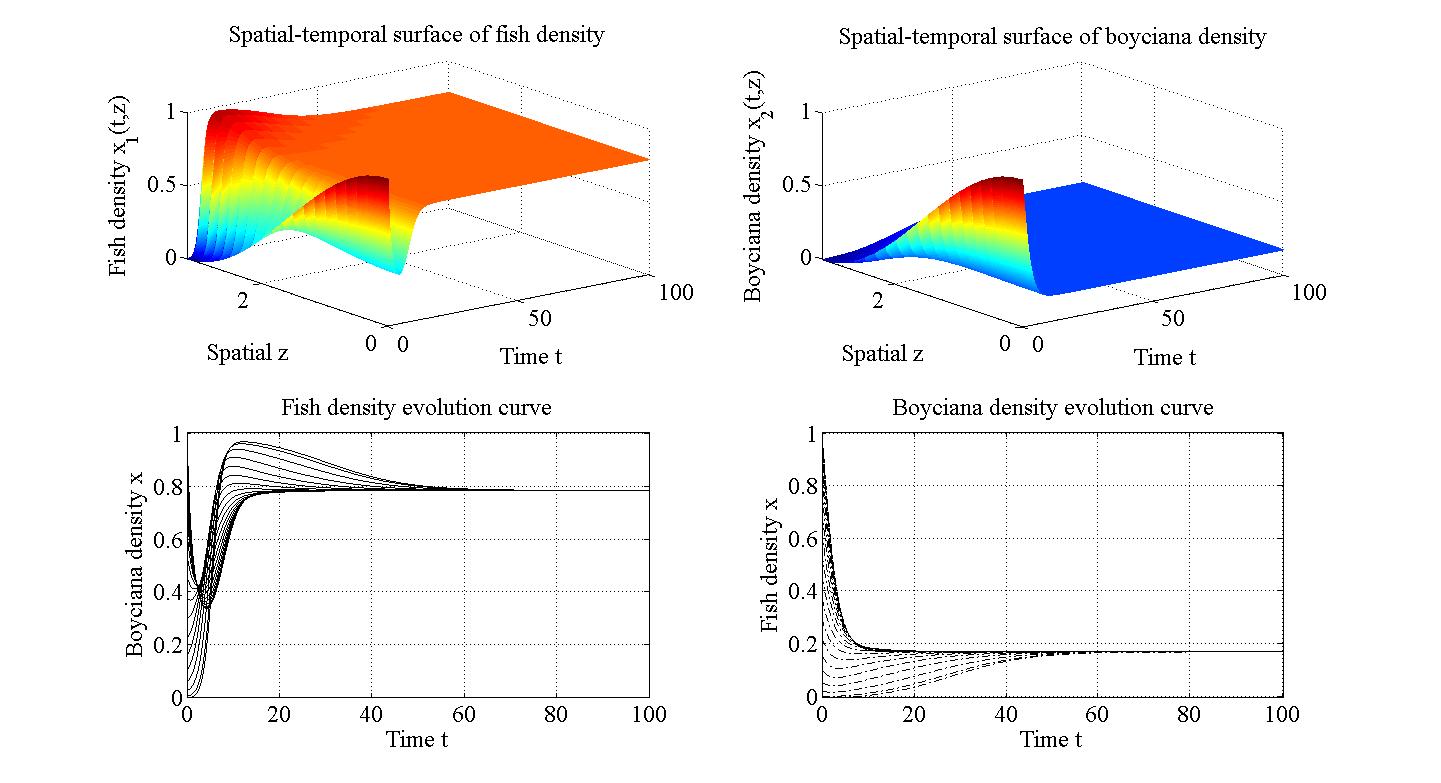}}
  \hspace{-0.05\textwidth}
  \subfigure[Spatial response curve]{
    \label{Fig_Sim_Stab_x3_subfig_b} 
    \includegraphics[width=0.5\textwidth]{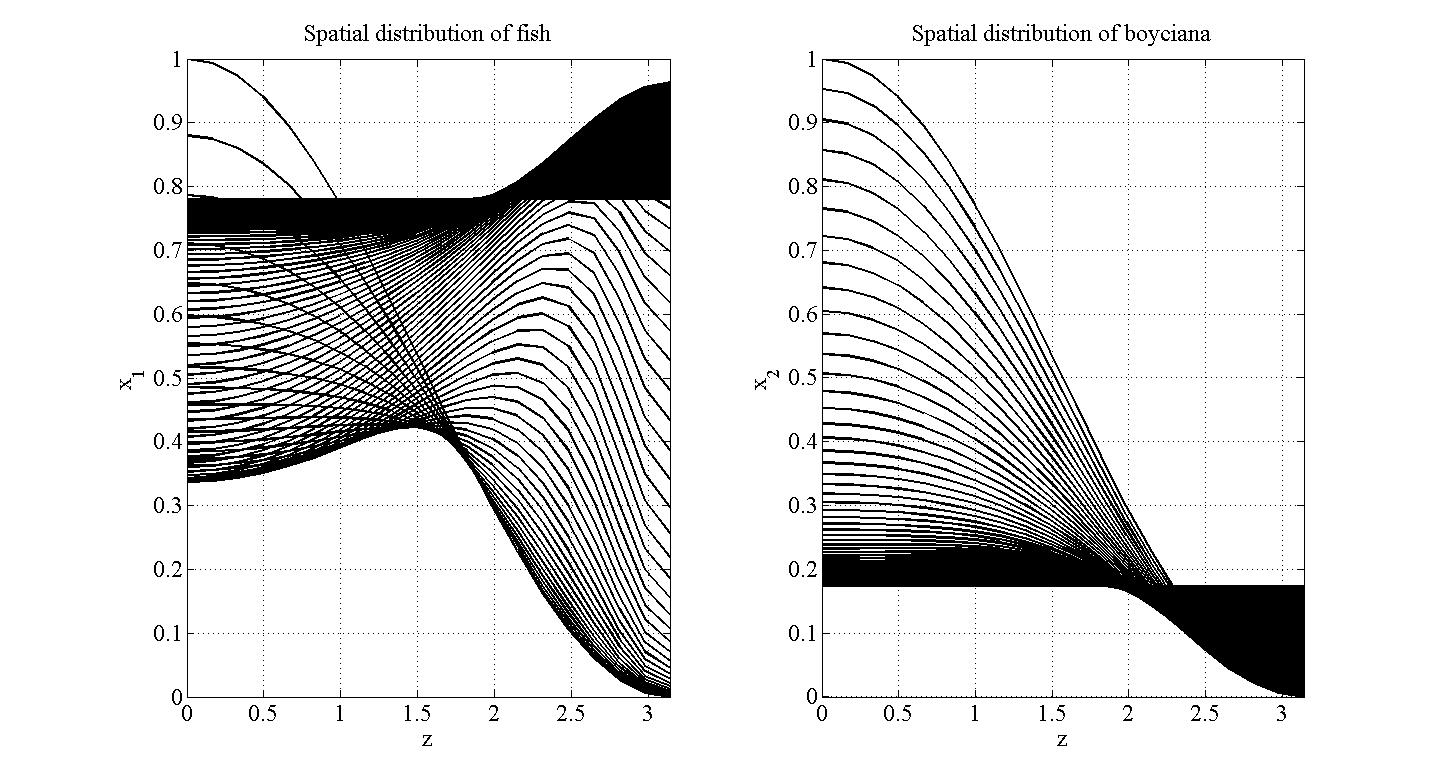}}
  \caption{Numerical results for the local stable property with $r=1.001$. (a) Stable steady state profiles and the corresponding time evolution curves. (b) Spatial response curve about $x_1(z)$ and $x_2(z)$ at differential discrete time points.}
  \label{Fig_Sim_Stab_x3} 
\end{figure}

\begin{figure}
  \centering
  \subfigure[Spatial-temporal profiles]{
    \label{Fig_Sim_Ustab_x3_subfig_a} 
    \includegraphics[width=0.5\textwidth]{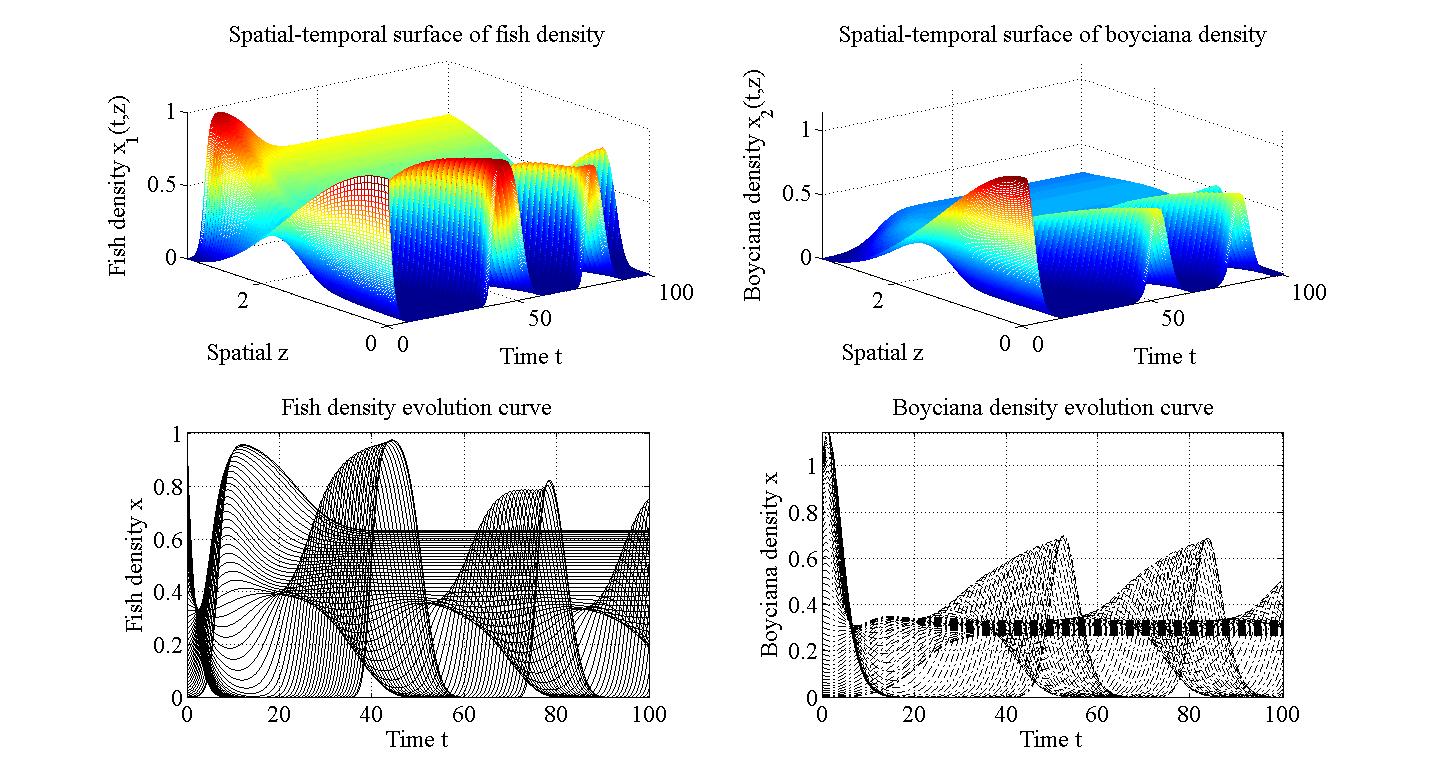}}
  \hspace{-0.05\textwidth}
  \subfigure[Spatial response curve]{
    \label{Fig_Sim_Ustab_x3_subfig_b} 
    \includegraphics[width=0.5\textwidth]{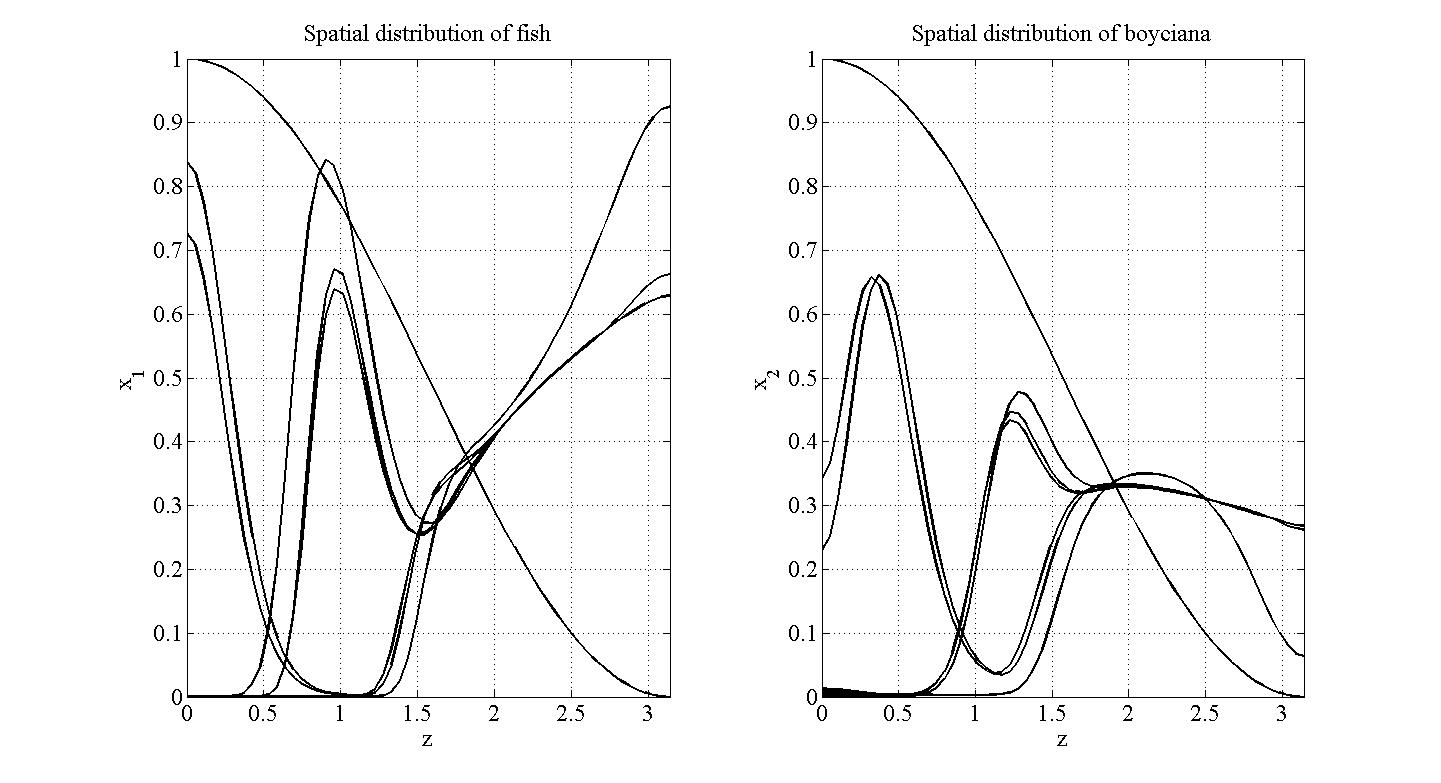}}
  \caption{Numerical results for the unstable property with human interference coefficient $r=100$. (a) Unstable steady state profiles and the corresponding time evolution curves. (b) Unstable spatial response curves about $x_1(z)$ and $x_2(z)$ at differential discrete time points.}
  \label{Fig_Sim_Ustab_x3} 
\end{figure}

In the system (\ref{sim_systems}), we fix all the system parameters except for $r$ as
\begin{eqnarray}\label{Sim_Parameter02}
d_1=d_2=0.01, c=1.000, \alpha=0.5000,\\
 d=0.3000, m=1.000, h_1=0.01, h_2=0.3.
\end{eqnarray}
As shown in Fig. \ref{Fig_Sim_Stab_x3} and Fig. \ref{Fig_Sim_Ustab_x3},  when  $r=1.001$ we can see that $x_1$ and $x_2$ are numerically converge to the desired steady state as $t\rightarrow\infty$, when $r=100$ $x_1, x_2$ can't maintain the steady state. This simulation result indicates that $r\rightarrow\infty$ may lead to the turning instable property of the PDAEs system.


\section*{Conclusion}
In this work, we have studied the dynamics of some boyciana-fish PDAEs system on wetland ecological system with human interference. Firstly, the existence and uniqueness properties of the positive solution of the elliptic human distribution model are investigated. Three solutions of this elliptic system are shown with difference practical significance. Moreover, the stability of the positive equilibrium and the persistence of the nontrivial positive solution of system (\ref{prey_predator systems}) are investigated by well-known linearized method and eigenvalue theory in PDEs. The 'diffusion-driven instability' property of the system is also discussed. It is shown that in a certain range of parameters, the positive constant steady state of (\ref{prey_predator systems}) is local asymptotically stable when the diffusion term $d_1, d_2$ satisfy certain condition and turning unstable when the conditions did not hold. What is more of interest is there exists some nonconstant steady state in our boyciana-fish model. And we propose some global stability analysis on boyciana and fish population by energy estimation. The theoretical energy estimation results show that the stability property of the boyciana-fish system is also determined by the parameter $r$. Here, $r$ is directly proportional to the human distribution quantity. Finally, we carry out some numerical results to illustrate the effectiveness of the development. And based on past fourteen years of bird data in Beidaihe wetland conservation district, by parameter estimation idea, we build maximum-minimum norm optimization algorithm to optimal the parameters of model (\ref{prey_predator systems}). With the practical data, our predict model is effective with $95.17\%$ average accuracy.
\section*{Conflict of Interests}
The authors declare that there is no conflict of interests
regarding the publication of this paper.
\section*{Acknowledgments}
The research is supported by N.N.S.F. of China
under Grant No. 61273008 and No. 61104003. The
research is also supported by the Key Laboratory of
Integrated Automation of Process Industry
(Northeastern University).

\bibliographystyle{unsrt}
\bibliography{mybib}
\end{document}